\documentclass[preprint, times]{elsarticle}
\usepackage{amsmath,amsfonts}
\usepackage{algorithmic}
\usepackage{algorithm}
\usepackage{array}
\usepackage[caption=false,font=normalsize,labelfont=sf,textfont=sf]{subfig}
\usepackage{textcomp}
\usepackage{stfloats}
\usepackage{url}
\usepackage{verbatim}
\usepackage{graphicx}

\usepackage{xcolor}
\usepackage{eurosym}
\usepackage{hyperref}
\usepackage[normalem]{ulem}
\usepackage{tikz}
\usepackage{tabularx}
\usetikzlibrary{positioning, calc}
\newdefinition{exmpl}{Example}

\begin{document}
	
	\title{State constrained stochastic optimal control of a PV system with battery storage via Fokker--Planck and Hamilton-Jacobi-Bellman equations}

	\author[1,2]{Alfredo Berm\'udez}
	\ead{alfredo.bermudez@usc.es}
	
	\author[1,2]{Iago Pad\'in\corref{cor1}}
	\ead{iago.padin.dominguez@usc.es}
	
	\cortext[cor1]{Corresponding author}
	
	\affiliation[1]{organization={Department of Applied Mathematics, University of Santiago de Compostela},
		city={Santiago de Compostela},
		country={Spain}}
	
	\affiliation[2]{organization={Galician Centre for Mathematical Research and Technology (CITMAga)},
		city={Santiago de Compostela},
		country={Spain}}
	
	\journal{arXiv}
	
	\begin{abstract}
		With the growing global emphasis on sustainability and the implementation of contemporary environmental policies, photovoltaic (PV) generation is playing an increasingly important role in modern power systems, while its intrinsic variability poses challenges for real-time operation and electricity market participation. This paper proposes a continuous-time stochastic optimal control framework for the joint optimization of real-time battery management and day-ahead market bidding of PV plants with energy storage. Solar irradiance, electricity prices, and battery dynamics are modeled through stochastic differential equations (SDEs), leading to a constrained stochastic control problem characterized by a coupled Hamilton-Jacobi-Bellman (HJB) and Fokker-Planck (FP) formulation. To mitigate the associated computational burden, a dimension-reduction strategy is introduced by decomposing the state space into controllable and uncontrollable components, yielding lower-dimensional optimality systems while preserving the continuous-time Dynamic Programming structure. Numerical results show that the reduced formulations achieve substantial computational savings, enabling real-time applicability without significant loss of performance. The proposed methodology is benchmarked against two rule-based strategies and a stochastic Model Predictive Control (MPC) approach, highlighting a favorable trade-off in terms of economic performance, computational efficiency, and suitability for day-ahead market participation.
	\end{abstract}
	
	\newpageafter{abstract}
	
	
	\maketitle
	
	\pagebreak
	
	\section*{Nomenclature}
	\addcontentsline{toc}{section}{Nomenclature}
	
	\subsection*{Notations}
	
	\noindent
	\begin{tabularx}{\linewidth}{lX}
		\hline
		\multicolumn{2}{l}{\textbf{General variables}} \\
		\hline
		$\mathcal{A}$ & Convex set defining state constraints \\
		$\mathbf{B}_s$ & Independent standard Brownian motions \\
		$\mathbf{b}$ & Drift field \\
		$f$ & Running-cost function in objective function \\
		$\mathbf{G}$ & State constraints function \\
		$g$ & Final time cost in objective function \\
		$J$ & Objective function \\
		$\mathcal{L}_{aug}$ & Augmented Lagrangian function \\
		$\mathrm{P}_\mathcal{A}$ & Orthogonal projection function on set $\mathcal{A}$ \\
		$s$ & Time (h) \\
		$t, T$ & Initial and final time instants (h) \\
		$\boldsymbol{\mathcal{U}}$ & Admissible control set \\
		$\mathbf{U}$ & Compact set for admissible controls \\
		$\boldsymbol{u}_s$ & Control \\
		$\boldsymbol{\check{u}}$ & Feedback control \\
		$\boldsymbol{u}(s)$ & Expectation of control \\
		$v$ & Value function \\
		$\mathbf{X}_s^{t,\mathbf{x}}$ & State of the system at time $s$, starting  from $\mathbf{x}$ at time $t$ \\
		$\mathbf{x}$ & Initial state \\ 
		$\mathbf{y}$ & Spatial variables \\
		$\zeta^*, \eta^*, \tau^*$ & Tolerances in algorithm \\
		$\lambda$ & Regularization parameter \\
		$\boldsymbol{\mu}(s)$ & Lagrange multiplier \\
		$\rho$ & Relaxation parameter \\
		$\mathbf{\Sigma}$ & Volatility tensor \\
		$\varphi$ & Probability density function \\
		$\varphi_\mathbf{x}$ & Initial condition for FP equation \\
		$\chi_{\mathcal{A}}$ & Indicator function of convex set $\mathcal{A}$ \\
		$\partial \chi_\mathcal{A}$ & Subdifferential of the indicator function \\
		$\omega_1, \omega_2$ & Decreasing parameters in algorithm \\
		\hline
	\end{tabularx}
	
	\noindent
	\begin{tabularx}{\linewidth}{lX}
		\hline
		\multicolumn{2}{l}{\textbf{PV--battery problem variables}} \\
		\hline
		$A_{PV}$ & Area of the PV system (m$^2$) \\
		$B_{Z,s}, B_{\Pi,s}, B_{\mathcal{E},s}$ & Standard Brownian motions \\
		$c_T$ & Terminal cost (\euro/MWh) \\
		$D$ & Day of the year (-) \\
		$\mathcal{E}_{max}, \mathcal{E}_{min}$ & Maximum and minimum battery capacity (MWh) \\
		$\mathcal{E}_s$ & Energy stored in the battery at time $s$ (MWh) \\
		$\mathcal{E}^0$ & Initial condition of the SDE for $\mathcal{E}_s$ (MWh) \\
		$I_s$ & Global Horizontal Irradiance at time $s$ (W/m$^2$) \\
		$I_{CS}(s)$ & Clear Sky Irradiance (W/m$^2$) \\
		$\hat{P}(I_s)$ & Deterministic function defining power produced by the PV plant (MW) \\
		$P_{\text{bat},s}$ & Power delivered by the battery at time $s$ (MW) \\
		$\check{P}_{\text{bat},s}$ & Power delivered by the battery at time $s$ in feedback form (MW) \\
		$P_{\text{grid},s}$ & Power sold to (positive) or purchased from (negative) the grid at time $s$ (MW) \\
		$\tilde{P}_{\text{grid}}(s)$ & Power bid at day-ahead electricity market auctions (MW) \\
		$P_{max}, P_{min}$ & Maximum (positive) and minimum (negative) power delivered by the battery (MW) \\
		$P_{\text{solar},s}$ & Power produced by the PV plant at time $s$ (MW) \\
		$Z_s$ & Clear Sky Index at time $s$ (-) \\
		$Z^0$ & Initial condition of the SDE for $Z_s$ (-)\\
		$\alpha_z(s)$ & Zenith solar angle (rad) \\
		$\beta$ & Tilt angle of the PV module from the horizontal (rad) \\
		$\gamma(s)$ & Solar inclination angle (rad) \\
		$\eta$ & Efficiency of the PV system (-) \\
		$\theta_Z(s)$ & Mean reverting function associated with $Z_s$ (-) \\
		$\theta_\Pi(s)$ & Mean reverting function associated with $\Pi_s$ (\euro/MWh) \\
		$\kappa_Z(s), \kappa_{\Pi}$ & Mean reverting parameters (h$^{-1}$) \\
		$\Pi_s$ & Electricity price at time $s$ (\euro/MWh) \\
		$\Pi^0$ & Initial condition of the SDE for $\Pi_s$ (\euro/MWh)\\
		$\Pi_{max}, \Pi_{min}$ & Price thresholds for rule-based method (\euro/MWh) \\
		$\sigma_{i,j}$ & Constant parameters for $i, j = Z, \Pi, \mathcal{E}$ (h$^{1/2}$) \\
		$\Psi$ & Non-linear function \\
		\hline
	\end{tabularx}
	
	\subsection*{Abbreviations}
	\noindent
	\begin{tabularx}{\linewidth}{lX}
		ADMM & Alternating Direction Method of Multipliers \\
		AR & Autoregressive \\
		ARMA & Autoregressive Moving Average \\
		CSI & Clear-Sky-Index \\
		DRO & Distributionally Robust Optimization \\
		DPP & Dynamic Programming Principle \\
		FP & Fokker--Planck \\
		GHI & Global Horizontal Irradiance \\
		HJB & Hamilton-Jacobi-Bellman \\
		LQ & Linear-quadratic \\
		MDP & Markov Decision Parameter \\
		MILP & Mixed-Integer Linear Programming \\
		MPC & Model Predictive Control \\
		PDE & Partial Differential Equation \\
		pdf & Probability density function \\
		PV & Photovoltaic \\
		RL & Reinforcement Learning \\
		SDE & Stochastic Differential Equation \\
		SPMP & Stochastic Pontryagin Maximum Principle \\
		TOU & Time-of-use \\
		VPP & Virtual Power Plant
	\end{tabularx}

	\section{Introduction}
	Renewable energies are facing unprecedent growth with the aim of achieving  the goal of zero net emissions by 2050 among other environmental objectives. In particular, solar photovoltaics (PV) are expected to be the leading renewable energy source in the future due to their scalability, declining costs, and minimal environmental impact. \cite{IEA1,IEA2}. Despite the evident advantages, this energy source also poses some technical challenges. The most critical inconvenient is the intermittency of solar irradiance. This issue complicates grid integration, energy planning, and market participation, making it necessary to develop control and optimization strategies to mitigate its effects.
	
	Advancements in grid interconnection and management, along with improvements in weather forecasting, are of paramount importance to overcome the challenges associated with renewable energy. In the case of PV systems, various strategies have been proposed for optimizing its integration \cite{SAA, TGSW}, including hybrid approaches combining PV panels with wind turbines \cite{MTIBG, BK}. More recently, the combination of PV plants with energy storage systems has arised as a promising alternative with effective results \cite{LLWWLZ}. In this context, batteries allow to smooth out sudden fluctuations in production, enhancing grid stability, improving overall energy reliability, and guaranteeing security of supply. Henceforth, this work focuses on a PV energy plant coupled with a battery system which participates in the electricity market.
	
	Accurately predicting the evolution of PV generation, energy market prices and battery dynamics requires precise mathematical modelling. The literature offers various approaches to this end: Markov decision processes (MDPs) are commonly used \cite{BPP}, especially as a baseline when working in reinforcement learning (RL) frameworks \cite{CKLZ, ZT}. Time series models, as autoregressive (AR) models \cite{CMSS} or autoregressive moving average (ARMA) models \cite{WPMLBLS}, are also applied in different methodologies for forecasting or decision-making. In this paper, stochastic differential equations (SDEs) are taken as the tool to handle uncertainty. These equations are widely used for the modelling of power systems, \cite{VAKB}. In the case of solar power generation, several models have been proposed, \cite{IMMM}. Here, the SDE introduced in \cite{BGGK} is adopted for the renewable factor, namely, a quotient of  solar irradiances. This equation relies on the concept of Clear Sky Irradiance, which was introduced in \cite{BH}. To model energy market prices, Johnson et al.~\cite{JHD} introduced a SDE with mean-reverting behavior, capturing the typical fluctuations and daily trends of electricity prices. For the evolution of the energy content of the battery, a bucket-type model is applied, which reflects the charge--discharge dynamics in a simplified yet effective manner, as described in~\cite{ZPPBAZBP}.
	
	Different microgrid scales and objectives were addressed in the literature within the context of PV systems operation and energy management optimization, thus leading to a broad range of mathematical frameworks for problem formulation and solving. A first intuitive control strategy consists of rule-based algorithms, which rely on predefined logical conditions to govern system operation without requiring complex optimization. Efstratiadi et al. \cite{ETPS} applied this technique to reduce the energy costs of residential users equipped with PV--battery systems. However, although computationally appealing, these approaches lack optimality guarantees and fail to provide feedback policies that adapt continuously to the system’s state.
	
	Alternatively, Model Predictive Control (MPC) approaches, provide a more efficient methodology to anticipate future conditions and optimize system operation over a prediction horizon while handling constraints explicitly \cite{ZPPBAZBP}. Wang et al. \cite{WNK} minimized the operation costs of a system comprising PV arrays, batteries and electric vehicles combining MPC with the alternating direction method of multipliers (ADMM) for handling constraints. Wu et al. \cite{WPMLBLS} developed a home energy management system model in the presence of PV, where control decisions were computed by combining MPC with a mixed-integer linear programming (MILP) framework. Despite their effectiveness in real-time operation, MPC-based approaches may be less suitable for day-ahead electricity market participation, where bidding decisions must be determined in advance and remain robust to future uncertainty. From the perspective of mathematical optimization, MILP is widely used in PV--battery energy management due to its ability to model discrete operational decisions, such as charging modes, or market participation, while optimizing system performance under linear constraints and objectives \cite{CMSS, QDRZSK}. Herding et al. \cite{HRJCP, HRJECP} integrated MILP with stochastic programming to control the uncertainties associated with PV--battery systems participation in day-ahead market auctions. 
	
	Other stochastic programming frameworks, such as those based on scenario generation or Monte Carlo sampling, can also be used in this setting and are commonly addressed using metaheuristic methods such as genetic algorithms \cite{ZY}. A related strand of control work employs stochastic programming with scenario trees, which discretizes the evolution of uncertain parameters into a finite number of scenarios organized in a tree structure. Although this representation provides a detailed view of uncertainty, it typically leads to rapidly increasing computational complexity as the number of scenarios grows, making real-time implementation challenging \cite{BCo}.
	
	Complementing traditional optimization and control techniques, recent research has developed methods designed to manage uncertainty in alternative ways. Instead of relying on a fully specified probability distribution of uncertain variables, Distributionally Robust Optimization (DRO) assumes that the true distribution belongs to an ambiguity set constructed from historical data or moment information, thus providing solutions that remain reliable under distributional shifts \cite{SYD}. Feng et al. \cite{FRWZ} applied DRO for the bidding strategy of a Virtual Power Plant (VPP) that participates in energy markets and handled constraints in the style of chance constraint programming framework introduced in \cite{CC}. Under both approaches, constraints are required to be satisfied at a prescribed probability level and under the worst-case probability distribution. More recently, RL has been applied to energy markets and renewable resource management, allowing the system to learn optimal control or bidding strategies from interactions with the environment \cite{ZT}. Chen et al. \cite{CKLZ} analyzed the operation of a grid-connected PV--battery system and demonstrated that Deep Reinforcement Learning methods outperform conventional rule-based control schemes in terms of operational efficiency.
	
	In this work, Dynamic Programming Principle (DPP) \cite{B} emerges as an appropriate framework for the goal of real-time operation and day-ahead market participation. Unlike some of the previously discussed approaches, which are predominantly data-driven (e.g., DRO), the methodology applied here is fully model-based, and in particular explicitly quantifies uncertainty through SDEs. Several studies have explored the application of DPP in energy systems from a discrete-time perspective demonstrating its effectiveness. Belloni et al. \cite{BPP} addressed the problem of energy management in a wind turbine facility with a battery under uncertainty. Zou et al. \cite{ZPLLYY} focused on a grid connected PV--battery system and demonstrated the better economic performance of the stochastic Dynamic Programming strategy compared with two different rule-based methods. Riffonneau et al. \cite{RBBP} paid special attention to peak shaving and took into account battery ageing for market interaction. Chouaf et al. \cite{CAA} applied DPP in a deterministic case proving its effectiveness when compared with a rule-based algorithm. However, to the best of authors' knowledge, no previously study has addressed the operation of PV plants with energy storage systems in a continuous-time setting. Modern electricity markets are increasingly shifting toward finer temporal resolutions. For example, the Iberian electricity market has recently moved from hourly to quarter-hourly trading intervals and now also includes a continuous market. In this context, continuous-time models provide a natural and flexible framework that remains consistent across different levels of market granularity, avoiding the need for repeated model reformulation.
	
	This paper covers this gap by employing stochastic optimal control theory within a continuous-time framework, which naturally leads to the Hamilton-Jacobi-Bellman (HJB) partial differential equation \cite{BCa}. For that purpose, a general problem formulation is proposed, making it valid from small to large PV energy producers and leaving aside the microgrid connections between PV plants and batteries or grid, which can be AC or DC. Nevertheless, it should be noted that the inclusion of state constraints such as bounds on battery energy levels introduces additional complexity, making the HJB equation alone insufficient for solving the problem. Fortunately, this difficulty can be addressed by reformulating the stochastic optimal control problem as a deterministic optimal control problem using the Fokker--Planck (FP) partial differential equation (PDE) \cite{ABNT, RWB}. This methodology has been successfully applied by the authors in previous work \cite{BP}, where a wind power plant with battery storage system was analyzed. This work extends this strategy to the case of a PV-storage system, introducing a key methodological difference in the treatment of state constraints. Whereas earlier formulations relied on a standard Lagrangian approach, we employ an Augmented Lagrangian formulation [37]. This modification is not merely technical: it strengthens the enforcement of state constraints and enables algorithmic improvements that enhance numerical stability and accelerate convergence of the optimality system.
	
	Besides covering the gap of applying continuous-time DPP to state constrained stochastical optimal control of PV plants with energy storage systems, this paper proposes an interesting novelty which lies in the reduction of the dimension of the stochastic control problem. Indeed, when the dimension of the problem grows, it suffers from the \textit{curse of dimensionality}, a concept introduced by Richard Bellman \cite{B2}. The proposed strategy mitigates the dimensionality challenge while maintaining consistency with the continuous-time Dynamic Programming Principle framework. Specifically, the state of the system can be decomposed into two components: one influenced by the control input and another evolving independently of it. Consequently, the FP equation \cite{Ri} can be decomposed into two equivalent lower-dimensional FP equations. This idea makes the approach highly scalable for practical applications.
	
	To further assess the practical relevance and performance of the proposed continuous-time Dynamic Programming framework, the methodology is benchmarked against representative alternatives from the literature. In particular, we consider two widely used rule-based strategies, namely a price-threshold policy and a time-of-use (TOU) policy, which are commonly adopted in practice due to their simplicity and low computational cost. In addition, we include a stochastic MPC approach, \cite{MKD, HPOM}, combined with the stochastic Pontryagin Maximum Principle (SPMP), \cite{P}, as a more advanced optimization-based benchmark, capable of incorporating uncertainty and constraints. This comparative analysis aims to highlight the trade-offs between control optimality, computational effort, and implementation complexity, and to position the proposed approach with respect to established methodologies used for energy management and market participation of renewable generation systems.
	
	Summarizing, this paper proposes a strategy to maximize the revenue of a renewable energy plant by optimizing both the real-time control of the storage system and the bids submitted to the preceding electricity market auction. By combining an Augmented Lagrangian treatment of state constraints with a principled dimension reduction strategy, the proposed methodology enhances both numerical efficiency and practical applicability. The remainder of the paper is organized as follows: Section~\ref{sect:II_stoch_problem} introduces the stochastic control problem to be addressed --both for a PV plant and in a general formulation that  serves as a baseline-- and derives its reformulation as a deterministic optimal control problem.  Section~\ref{sect:III_opt_system} describes a method for handling state constraints and an algorithm for solving the resulting optimality system. Section~\ref{sect:IV_space_red} discusses spatial dimension reduction. Section~\ref{sect:V_results} returns to the PV case study and reports numerical results, benchmarking the proposed methodology against two representative rule-based approaches and a stochastic MPC framework. Finally, Section~\ref{sect:VI_conclusion} concludes the study.

	\section{Stochastic control problem} \label{sect:II_stoch_problem}
	\noindent This section introduces a stochastic control problem designed to handle the challenges of a PV energy plant with energy storage systems. The problem is then generalized to a broader formulation and rewritten as a deterministic control problem.

	\subsection{Mathematical modelling} \label{sect:II_I_PV_problem}
\begin{figure}[!t]
		\centering
		\includegraphics[width=3.in]{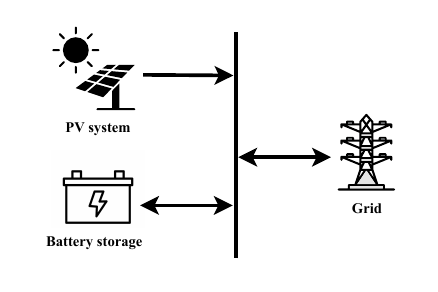}
		\caption{Scheme of system considered and feasible power flows.}
		\label{fig:scheme}
	\end{figure}

	\noindent The model of the system consists of three stochastic processes which are the Clear-Sky-Index (CSI), $Z_s$, the price of electricity, $\Pi_s$, and the energy content of the battery, $\mathcal{E}_s$. Fig.~\ref{fig:scheme} illustrates a schematic representation of the assumed connections between the system components and the electrical grid.
	
	The CSI was introduced in \cite{BGGK} to characterize the renewable factor associated with solar irradiance. It is defined as the ratio between the Global Horizontal Irradiance (GHI), $I_s$, and the Clear Sky Irradiance, $I_{CS}$, which represents the maximum irradiance observable under cloud-free conditions:
	\begin{equation}
		\label{eq:Z_s}
		Z_s = \frac{I_s}{I_{CS}(s)}.
	\end{equation}
	Notice that $I_s$, and consequently $Z_s$, are stochastic processes, while the deterministic function $I_{CS}$ can be computed using the empirical formula introduced in \cite{BGGK}, which is given for a specific day of the year $D$ and a zenith solar angle $\alpha_z(s)$:
	\begin{equation}
		\label{eq:I_CS}
		I_{CS}(s) = \cos(\alpha_z(s))^{1.2} \bigg( 83.69 \sin\left(\frac{2\pi(D + 82.07)}{365.24}\right)
		+ 1130.44 \bigg).
	\end{equation}
	
	To model the state vector, a stochastic system is proposed. For the CSI, the mean reverting SDE introduced in \cite{BGGK} is used, while the electricity price follows a mean-reverting-like SDE as in \cite{JHD}. These models capture random fluctuations while, on average, reflecting a natural tendency to revert to equilibrium. Finally, a bucket model is adopted to represent the energy stored in the battery, describing the charging--discharging dynamics as inflows and outflows of energy. Thus, the state equations are
	\begin{subequations} \label{eq:SDE_PV}
		\begin{align}
			\label{eq:SDE_Z}
			dZ_{s} 	=&\kappa_{Z}\left(\theta_{Z}(s)
			-Z_{s}\right)ds + \sigma_{Z,Z}Z_s(1-Z_s)dB_{Z,s},
			\\
			\nonumber
			d\Pi _s =&\kappa_{\Pi}\left(\theta_{\Pi}(s)+\frac{1}{\kappa_{\Pi}}\frac{d\theta_\Pi}{ds}(s)-\Pi _s\right)ds 
			\\
			\label{eq:SDE_Pi}
			+& \sigma_{\Pi,Z}\Pi_sdB_{Z,s} + \sigma_{\Pi,\Pi}\Pi_s dB_{\Pi,s},
			\\
			\label{eq:SDE_E}
			d\mathcal{E}_{s} =&-P_{\text{bat},s}ds + \sigma_{\mathcal{E},\mathcal{E}}\mathcal{E}_{s} dB_{\mathcal{E},s},
			\\
			\label{eq:SDE_ci}
			Z_{0} =& Z^0, \quad \Pi _0= \Pi^0, \quad \mathcal{E}_{0}=\mathcal{E}^0.
		\end{align}
	\end{subequations}
	This model includes the following parameters and functions:
	\begin{itemize}
		\item Mean reverting functions $\theta_Z$ (-) and $\theta_\Pi$ (\euro/MWh).
		\item Mean reverting parameters $\kappa_{Z}$ (h$^{-1}$) and $\kappa_{\Pi}$ (h$^{-1}$).
		\item Standard Brownian motions $B_{Z,s}$, $B_{\Pi,s}$ and $B_{\mathcal{E},s}$.
		\item Constant parameters $\sigma_{i,j}$, $i,j = Z, \Pi, \mathcal{E}$ (h$^{1/2}$).
	\end{itemize}
	
	In the mathematical problem, $P_{\text{bat},s}$ is taken as the control function with the sign convention: $P_{\text{bat},s}$ is positive when the battery is discharging and negative when it is charging. 
	
	The objective of the stochastic control problem is to maximize the expected revenues of the system. It corresponds to the expected value of the electricity price times the net energy exchanged with the grid over the interval of interest, namely,
	\begin{equation}
		\label{eq:Fun_Obj}
		J(P_{\text{bat},s})=\mathrm{E}\Bigg[\int_0^T \Pi_s P_{\text{grid},s}\,\mathrm{d}s + g(\mathcal{E}_{T})\Bigg].
	\end{equation}
	Here, $g: \mathbb{R} \rightarrow \mathbb{R}$ is a deterministic function representing the cost of the battery's energy content at the end of the process. This term ensures that the battery is not completely depleted at the final time, preserving some energy for the start of the next operational period.
	
	The power exchanged with the grid is the sum of the power delivered by the battery and the power produced by the PV system, i.e., 
	$$P_{\text{grid},s} = P_{\text{solar},s} + P_{\text{bat},s} = \hat{P}(I_s)+~P_{\text{bat},s},$$
	where the power produced by the PV plant can be computed by the deterministic function
	\begin{equation}
		\label{eq:PV_power}
		\hat{P}(I)= 10^{-6} A_{PV} \eta I \frac{\sin(\gamma(s) + \beta)}{\sin \gamma(s)}.
	\end{equation}
	This function depends on the area of the PV system $A_{\text{PV}}$ (m$^2$), the efficiency $\eta$ (-), the tilt angle of the PV module from the horizontal $\beta$ (rad) and the solar inclination angle~$\gamma(s)$~(rad).
	
	Finally, some constraints must be considered. On the one hand, there are limitations on the power delivered by the battery: 
	$$P_{min}\leq P_{\text{bat},s} \leq P_{max}.$$ 
	On the other hand, the state variable 
	$\mathcal{E}_s$ is subject to physical constraints, as it must remain within the minimum and maximum capacity limits of the battery, i.e, $$\mathcal{E}_{min} \leq \mathcal{E}_s \leq \mathcal{E}_{max}.$$ 
	This state constraint is reformulated as an expectation, that is, 
	$$\mathcal{E}_{min}\leq \mathrm{E}[\mathcal{E}_{s}]\leq \mathcal{E}_{max}.$$

	\subsection{General framework}
	\noindent The above problem belongs to the following general stochastic control framework. First,
	the state of the system is a vector stochastic process, $\mathbf{X}_s^{t,\mathbf{x}}\in\mathbb{R}^n$, with $s \in [t, T]$ which satisfies the SDE:
	\begin{subequations}\label{eq:SDE_system}
		\begin{align}
			\label{eq:SDE_gral}
			d\mathbf{X}^{t,\mathbf{x}}_s &=\mathbf{b}(s,\mathbf{X}^{t,\mathbf{x}}_s, \boldsymbol{u}_s)ds +\mathbf{\Sigma}(s,\mathbf{X}^{t,\mathbf{x}}_s)d\mathbf{B}_s,
			\\
			\label{eq:SDE_gral_IC}
			\mathbf{X}^{t,\mathbf{x}}_t&=\mathbf{x}.
		\end{align}
	\end{subequations}
	Here, $\boldsymbol{u}_s$ denotes the control of the problem and verifies that $\boldsymbol{u} \in \boldsymbol{\mathcal{U}}$ where 
	$$\boldsymbol{\mathcal{U}} = \{\boldsymbol{u}:[t,T]\rightarrow \mathrm{\mathbf{U}} \subset \mathbb{R}^r| ~ \boldsymbol{u} ~ \text{measurable}\},$$
	being $\mathrm{\mathbf{U}}$ a compact set. In the above SDE, $\mathbf{b}$ and $\mathbf{\Sigma}$ are given vector and tensor fields, respectively, while the $d$ components of $\mathbf{B}_s$ are independent standard Brownian motions.
	
	The objective function, also known as cost-to-go, is defined in terms of a running cost, $f$, accumulated over the time horizon, and a terminal cost, $g$, that accounts for the final state:
	\begin{equation}
		\label{eq:Fun_Obj_Stoch_gral}
		J(t,\mathbf{x},\boldsymbol{u})= \mathrm{E}\left[\int_t^T f(s,\mathbf{X}_s^{t,\mathbf{x}},\boldsymbol{u}_s)\,\mathrm{d}s+ g(\mathbf{X}_T^{t,\mathbf{x}})\right].
	\end{equation}
	Finally, let $\mathcal{A}$ be a convex subset of $\mathbb{R}^n$, then the state constraints are formulated as 
	$$\mathrm{E}[\mathbf{X}_s^{t,\mathbf{x}}] \in \mathcal{A}, \, s \in [t,T].$$

	\subsection{Stochastic control problem reformulation}
	\noindent It is well-known that, in the absence of state constraints, the preceding stochastic optimal control problem can be solved using the DPP which is equivalent to the HJB equation. The equation encodes the principle of optimality, connecting the system's dynamics with the \textit{value function} --the maximal performance achievable starting from any time and state. However, to handle state constraints it is convenient to reformulate the original problem as a deterministic optimal control problem with the FP equation playing the role of (deterministic) state equation. This procedure and following steps are summarized in Figure \ref{fig:diagram} and further detailed in the subsequent sections.

	\begin{figure}[!t]
		\centering
		\begin{tikzpicture}[scale=0.8, transform shape, node distance=1cm and -2cm, auto]
			\node (stoch pb) [rectangle, draw, rounded corners, align=center] {State constrained stochastic optimal control problem};
			\node (fp) [rectangle, draw, rounded corners, align=center, below left=of stoch pb] {FP equation \\ Probability density function};
			\node (hjb) [rectangle, draw, rounded corners, align=center, below right=of stoch pb] {HJB equation \\ Feedback control};
			
			\coordinate (midpoint) at ($(fp)!0.5!(hjb)$);
			
			\node (deter pb) [rectangle, draw, rounded corners, align=center, below=of midpoint] {State constrained deterministic optimal control problem};
			\node (Laug) [rectangle, draw, rounded corners, align=center, below=of deter pb] {Augmented Lagrangian approach \\ Lagrange multiplier};
			\node (opt) [rectangle, draw, rounded corners, align=center, below=of Laug] {Optimality system \\ FP, HJB  and Lagrange multiplier characterization};
			
			\draw[->, thick, >=stealth] (stoch pb.south) -- (fp.north);
			\draw[->, thick, >=stealth] (stoch pb.south) -- (hjb.north);
			\draw[->, thick, >=stealth] (fp.south) -- (deter pb.north);
			\draw[->, thick, >=stealth] (hjb.south) -- (deter pb.north);
			\draw[->, thick, >=stealth] (deter pb.south) -- (Laug.north);
			\draw[->, thick, >=stealth] (Laug.south) -- (opt.north);
			
		\end{tikzpicture}
		\caption{Workflow of the stochastic optimal control methodology with state constraints.}
		\label{fig:diagram}
	\end{figure}

	First, the FP equation (also known as the forward Kolmogorov equation) is a linear PDE whose solution is the joint pdf of the stochastic process $\mathbf{X}_s^{t,\mathbf{x}}$ satisfying \eqref{eq:SDE_gral}, \eqref{eq:SDE_gral_IC}. Denoting the pdf of $\mathbf{X}_s^{t,\mathbf{x}}$ by 
	$$\varphi(s,.): \mathbf{y} \in \mathbb{R}^n \rightarrow\varphi(s,\mathbf{y})\in\mathbb{R},$$ 
	the function $\varphi$ satisfies the following initial-value problem:
	\begin{subequations} \label{eq:FP}
		\begin{align}
			\nonumber
			&\frac{\partial \varphi}{\partial s} (s,\mathbf{y})+\mbox{div}\Big(\mathbf{b}(s,\mathbf{y},\boldsymbol{u})\varphi(s,\mathbf{y})\Big)
			\\
			\label{eq:FP_eq}
			&\hspace{1.2cm}-\frac{1}{2}\mbox{div}\mathbf{div}\Big(\mathbf{\Sigma} (s,\mathbf{y})\mathbf{\Sigma} (s,\mathbf{y})^\mathtt{T}\varphi(s,\mathbf{y})\Big)=0,
			\\
			\label{eq:FP_IC}
			&\varphi(t,\mathbf{y})=\varphi_\mathbf{x}(\mathbf{y}),\ \mathbf{y}\in\mathbb{R}^n.
		\end{align}
	\end{subequations}
	From the SDE~\eqref{eq:SDE_system} one obtains this PDE, which characterizes the evolution of the state's probability density under uncertainty and applied controls. This perspective is key to stating constraints and achieving realistic control behavior.
	
	Second, we recall that the HJB equation provides the optimal control in feedback form. In other words, by solving this equation, we obtain in advance a deterministic feedback function that prescribes the optimal control based on the current state of the system, more specifically,
	$$\boldsymbol{\check{u}}:(s, \mathbf{y})\in  [t,T] \times \mathbb{R}^n \rightarrow \boldsymbol{\check{u}}(s, \mathbf{y}) \in \mathbb{R}^r,$$ 
	such that the optimal control at time $s$ is 
	$$\boldsymbol{u}_s = \boldsymbol{\check{u}}(s, \mathbf{X}_s^{t,\mathbf{x}}).$$ 
	Notice that, as $\boldsymbol{\check{u}}$ is a function of the state of the system, then $\boldsymbol{u}_s$ is a stochastic process.
	
	This approach is highly convenient for our objectives. Indeed, from the perspective of energy management, once the optimal feedback control $\boldsymbol{\check{u}}$ is computed, real-time implementation requires only measuring the current state of the system. On the other hand, if the optimal control is to be used in advance --as in the case of bidding in the day-ahead market-- a plausible solution is to use its expected value, i.e, 
	\begin{equation}
		\label{eq:exp_u}
		\boldsymbol{u}(s) = \mathrm{E}[\boldsymbol{u}_s] = \int_{\mathbb{R}
			^n} \boldsymbol{\check{u}}(s, \mathbf{y}) \varphi(s, \mathbf{y}) \, \mathrm{d}\mathbf{y}.
	\end{equation}
	
	\subsection{Deterministic control problem}
	\noindent Now we can reformulate the original stochastic problem with state constraints as a deterministic control problem. To do so, we consider the pdf $\varphi(s, \mathbf{y})$ as the state of the system. Accordingly, the state equation becomes the FP equation given by \eqref{eq:FP}. Furthermore, we take the deterministic feedback function $\boldsymbol{\check{u}}$ as the control. 
	
	Regarding the objective function to be maximized, the expectation involved can be computed using $\varphi(s, \mathbf{y})$. Therefore, it is rewritten as
	\begin{align}
		\nonumber
		J(\boldsymbol{\check{u}})&= \int_t^T \int_{\mathbb{R}^n} f(s,\mathbf{y},\boldsymbol{\check{u}}(s, \mathbf{y})) \varphi(s, \mathbf{y}) \, \mathrm{d}\mathbf{y}\mathrm{d}s
		\\
		\label{eq:Fun_Obj_Deter_gral}
		&+ \int_{\mathbb{R}^n} g(\mathbf{y})\varphi(T, \mathbf{y}) \, \mathrm{d}\mathbf{y}. 
	\end{align}
	Similarly, state constraints can also be rewritten using $\varphi(s, \mathbf{y})$:
	\begin{equation}
		\label{eq:mean_constr_deter_gral}
		\int_{\mathbb{R}^n}\mathbf{y} \varphi(s, \mathbf{y}) \, \mathrm{d}\mathbf{y} \in \mathcal{A}.
	\end{equation}
	
	According to Figure~\ref{fig:diagram}, the next step is to handle these state constraints via an Augmented Lagrangian approach and to derive the optimality system for the resulting state constrained deterministic optimal control problem. The following section details these steps.

	\section{Optimality system} \label{sect:III_opt_system}
	\noindent This section presents the optimality conditions for the previous deterministic optimal control problem. Some insights into their derivation are provided. For further details, the reader is referred to \cite{BP} where the formal derivation is carried out. It is worth noting that, while in \cite{BP} a standard Lagrangian function is used to handle the state constraints, the Augmented Lagrangian as described, for instance, in \cite{FG} is adopted here.
	
	The Augmented Lagrangian method adds a penalty term to the classical Lagrangian, improving conditioning and accelerating convergence in constrained problems. It typically delivers more robust numerics and avoids the excessively large penalties that can degrade accuracy in standard penalty methods. While the derivation follows a similar structure, the resulting formulation benefits from these advantages, making it better suited to our application.
	
	By treating the state equation as an additional constraint, the corresponding Lagrange multiplier will be denoted by $v(s, \mathbf{y})$. This multiplier turns to be the \textit{adjoint state} in the \textit{maximum principle} for deterministic optimal control problem. The adjoint state can be interpreted as a backward-propagating variable that quantifies the sensitivity of the final objective to changes in the trajectory of the system. While the system state evolves forward in time, the adjoint state evolves backward, establishing the connection between system behavior and optimality conditions. Regarding the state constraints, the notation $\boldsymbol{\mu}(s)$ will be used for the Lagrange multiplier associated with~\eqref{eq:mean_constr_deter_gral}. Thus, we are led to define the Augmented Lagrangian function:
	\begin{align}
		\nonumber
		\mathcal{L}&_{aug}(\boldsymbol{\check{u}},\varphi,v,\boldsymbol{\mu}) = \int_t^T\int_{\mathbb{R}^n}f(s,\mathbf{y},\boldsymbol{\check{u}}(s,\mathbf{y}))\varphi(s,\mathbf{y})\,\mathrm{d}\mathbf{y} \mathrm{d}s 
		\\
		\nonumber
		&+ \int_{\mathbb{R}^n}g(\mathbf{y})\varphi(T,\mathbf{y})\,\mathrm{d}\mathbf{y} - \int_{\mathbb{R}^n}v(T,\mathbf{y})\varphi(T,\mathbf{y}) \, \mathrm{d}\mathbf{y} 
		\\
		\nonumber
		&+\int_{\mathbb{R}^n}v(t,\mathbf{y})\varphi(t,\mathbf{y}) \,\mathrm{d}\mathbf{y}  +\int_t^T \int_{\mathbb{R}^n} \frac{\partial v}{\partial s} (s,\mathbf{y})\varphi(s,\mathbf{y})\mathrm{d}\mathbf{y} \,\mathrm{d}s
		\\
		\nonumber
		&+ \int_t^T  \int_{\mathbb{R}^n} \varphi(s,\mathbf{y}) \mathbf{b}(s,\mathbf{y},\boldsymbol{\check{u}}(s,\mathbf{y}))\cdot \mathrm{grad}\, v (s,\mathbf{y})\,\mathrm{d}\mathbf{y}\,\mathrm{d}s
		\\
		\nonumber
		&- \int_t^T\int_{\mathbb{R}^n}\frac{1}{2}\mathbf{div}\Big( \mathbf{\Sigma} (s,\mathbf{y})\mathbf{\Sigma} (s,\mathbf{y})^\mathtt{T} \varphi (s,\mathbf{y})\Big)\cdot \mathrm{grad}\, v (s,\mathbf{y})\,\mathrm{d}\mathbf{y}\,\mathrm{d} s
		\\
		\nonumber
		&- \int_t^T (\chi_{\mathcal{A}})_\lambda \left(\int_{\mathbb{R}^n}\mathbf{y}\varphi (s,\mathbf{y})\,\mathrm{d}\mathbf{y} + \lambda \boldsymbol{\mu}(s)\right)\mathrm{d}s 
		\\
		\label{eq:laug}
		&+ \frac{\lambda}{2} \int_t^T || \boldsymbol{\mu}(s)||^2 \mathrm{d}s.
	\end{align}
	
	The last two terms enforce the state constraint via an Augmented Lagrangian framework, using the indicator function of $\mathcal{A}$, $\chi_{\mathcal{A}}$, to ensure feasibility. Because $\chi_{\mathcal{A}}$ is non-differentiable, we  replace it with its Moreau-Yosida regularization, $(\chi_{\mathcal{A}})_\lambda$, where $\lambda>0$ is a regularization parameter (see \cite{BM} for further details).
	
	As shown in \cite{BP}, taking first-order partial derivatives of \eqref{eq:laug} yields the necessary optimality conditions for the deterministic optimal control problem, leading to the following optimality system
	{\allowdisplaybreaks
		\begin{subequations}\label{eq:opt_system}
			\begin{align}
				\nonumber
				&\frac{\partial \varphi}{\partial s} (s,\mathbf{y})+\mathrm{div}\Big(\mathbf{b}(s,\mathbf{y},\boldsymbol{\check{u}}(s,\mathbf{y}))\varphi(s,\mathbf{y})\Big) 
				\\*
				\label{eq:os11_1_laug}
				& \hspace*{1.15cm} -\frac{1}{2}\mathrm{div}\mathbf{div}\Big(\mathbf{\Sigma} (s,\mathbf{y})\mathbf{\Sigma} (s,\mathbf{y})^\mathtt{T}\varphi(s,\mathbf{y})\Big)=0,
				\\*
				\label{eq:os12_1_laug}
				&\varphi(t,\mathbf{y})=\varphi_\mathbf{x}(\mathbf{y}),\ \mathbf{y}\in\mathbb{R}^n,
				\\
				\nonumber
				&\frac{\partial v}{\partial s}(s,\mathbf{y}) +\max_{\boldsymbol{u}\in\mathrm{\mathbf{U}}}\left\{\mathbf{b}(s,\mathbf{y},\boldsymbol{u})\cdot\mathrm{grad}\, v(s,\mathbf{y})+f(s,\mathbf{y},\boldsymbol{u})\right\}
				\\*
				\nonumber
				&\hspace*{1.1cm}+\frac{1}{2}\mathbf{\Sigma} (s,\mathbf{y})\mathbf{\Sigma} (s,\mathbf{y})^\mathtt{T}\cdot \mathrm{Hess}\,(v)(s,\mathbf{y}) 
				\\*
				\label{eq:os13_1_laug}
				&\hspace*{1.1cm}- (\partial \chi_{\mathcal{A}})_\lambda \left(\int_{\mathbb{R}^n}\mathbf{y}\varphi (s,\mathbf{y})\,\mathrm{d}\mathbf{y} + \lambda \boldsymbol{\mu}(s)\right) \cdot \mathbf{y} = 0,
				\\*
				\label{eq:os14_1_laug}
				&  v(T,\mathbf{y})=g(\mathbf{y}),\ \mathbf{y}\in\mathbb{R}^n,
				\\
				\label{eq:os15_1_laug}
				&\boldsymbol{\mu}(s) = (\partial \chi_{\mathcal{A}})_\lambda \left(\int_{\mathbb{R}^n}\mathbf{y}\varphi (s,\mathbf{y})\,\mathrm{d}\mathbf{y} + \lambda \boldsymbol{\mu}(s)\right),
				\\
				\label{eq:os16_1_laug}
				& \boldsymbol{\check{u}}(s,\mathbf{y})=\underset{\boldsymbol{u}\in\mathrm{\mathbf{U}}}{\text{arg max}}\big\{\mathbf{b}(s,\mathbf{y},\boldsymbol{u})\cdot\mathrm{grad}\, v(s,\mathbf{y}) +f(s,\mathbf{y},\boldsymbol{u})\big\}.
			\end{align}
		\end{subequations}
	}
	
	The first optimality conditions \eqref{eq:os11_1_laug}, \eqref{eq:os12_1_laug} are the FP equation. The second optimality conditions correspond to the adjoint state equation, \eqref{eq:os13_1_laug}, \eqref{eq:os14_1_laug}, which is a backward nonlinear PDE. However, we emphasize that, for each given control $\boldsymbol{u}$, it is linear in $v$. Importantly, in the absence of state constraints, i.e., when $\boldsymbol{\mu}(s) = \mathbf{0}$, this equation reduces to the standard HJB equation for the unconstrained stochastic control problem and its solution, $v$, corresponds to the corresponding value function. Moreover, in this case, the HJB equation can be solved independently of the FP equation, which is no longer required, thus recovering the standard DPP.
	
	Regarding the third optimality condition \eqref{eq:os15_1_laug}, it characterizes the Lagrange multiplier $\boldsymbol{\mu}(s)$. Here, $\partial \chi_\mathcal{A}$ denotes the subdifferential of the indicator function of a set $\mathcal{A}$, \cite{Ro}. Its Moreau-Yosida regularization is 
	$$(\partial \chi_\mathcal{A})_\lambda = \frac{1}{\lambda}(\mathbf{y} - \mathrm{P}_\mathcal{A}(\mathbf{y}) ),$$ 
	being $\mathrm{P}_\mathcal{A}$ the orthogonal projection on the set $\mathcal{A}$. 
	Finally, the feedback control $\boldsymbol{\check{u}}$ is characterized by~\eqref{eq:os16_1_laug}.

	To build intuition for the methodology developed so far, we consider a simple one-dimensional stochastic control problem.
	\begin{exmpl}		
		Let the scalar stochastic state evolve as
		\begin{align*}
			dX_s^{t,x} &= u_s ds + \sigma dB_s, \\
			X_t^{t,x} &= x,
		\end{align*}
		where $u_s$ is the control input and $\sigma>0$ denotes the noise intensity. The objective is to minimize the expected cost given by
		\begin{equation*}
			J(t,x,u) = \mathrm{E}\left[\int_t^T \left(u_s^2 + (X_s^{t,x})^2 \right) ds + (X_T^{t,x})^2\right].
		\end{equation*}
		
		This formulation is the classical linear-quadratic (LQ) problem. In the unconstrained case,
		the optimal feedback control follows in closed form from the Riccati equation. Alternatively,
		the DPP yields a closed-loop solution via the
		HJB equation, which can be solved independently and takes the form
		\begin{align*}
			&\frac{\partial v}{\partial s}(s,y) +\max_{u}\left\{u \frac{\partial v}{\partial y} (s,y)+ u^2 + y^2\right\} + \frac{1}{2}\sigma^2 \frac{\partial^2 v}{\partial y^2} (s,y) = 0,
			\\
			&  v(T,y)=y^2,\ y\in\mathbb{R}.
		\end{align*}
		
		To illustrate the need for the proposed methodology, add simple bounds on the state expectation:
		\begin{equation*}
			\mathrm{E}[X_s^{t,x}] \in \mathcal{A} = [X_{min}, X_{max}], \quad s \in [t,T].
		\end{equation*}
		
		This prompts us to clarify the roles of the FP and HJB equations. First, the FP equation governs the evolution of the pdf \(\varphi\) of \(X_s^{t,x}\), enabling a deterministic reformulation of the problem:
		\begin{align*}
			&\frac{\partial \varphi}{\partial s} (s,y)+\frac{\partial}{\partial y} \left(u \varphi(s,y)\right) - \frac{1}{2}\frac{\partial^2}{\partial y^2}\left(\sigma^2 \varphi(s,y)\right)=0,
			\\
			&\varphi(t,y)=\varphi_x(y),\ y\in\mathbb{R}.
		\end{align*}
		Second, the HJB equation characterizes the value function \(v(s,y)\) and provides the feedback
		law \(u_s=\check u\big(s,X_s^{t,x}\big)\), which will be used as the deterministic control.
		
		With these tools, the stochastic LQ problem can be recast as a deterministic optimal control problem: the state is the pdf \(\varphi\), the state equation is the FP equation driven by the deterministic feedback \(\check u\), and the objective functional becomes
		\begin{equation*}
			J(\check{u}) = \int_\mathbb{R} \int_t^T \left(\check{u}^2(s,y) + y^2 \right) \varphi(s,y) \, ds \, dy + \int_\mathbb{R} y^2 \varphi (T,y) \, dy,
		\end{equation*}
		while the state bounds read
		$$ \int_{\mathbb{R}} y \varphi(s,y) \, d y \in \mathcal{A} = [X_{min}, X_{max}].$$
		
		To enforce the state constraints, we introduce a scalar Lagrange multiplier \(\mu(s)\). The resulting optimality conditions require solving not only an HJB-type equation but also the FP equation and the multiplier characterization. Specifically, in the constrained case the HJB equation reads
		\begin{align*}
			&\frac{\partial v}{\partial s}(s,y) +\max_{u}\left\{u \frac{\partial v}{\partial y} (s,y)+ u^2 + y^2\right\} + \frac{1}{2}\sigma^2 \frac{\partial^2 v}{\partial y^2} (s,y) \\
			&\hspace*{1.1cm}+ (\partial \chi_\mathcal{A})_\lambda \left(\int_{\mathbb{R}} y \varphi(s, y) \, dy + \lambda \mu(s)\right) y = 0,
			\\
			&  v(T,y)=y^2,\ y\in\mathbb{R},
		\end{align*}
		where the last term is nonnegative in the minimization setting. The Lagrange multiplier satisfies
		$$ \mu(s) = (\partial \chi_\mathcal{A})_\lambda \left(\int_{\mathbb{R}} y \varphi(s, y) \, dy + \lambda \mu(s)\right).$$
		
		Although highly simplified, this constrained LQ example highlights the joint roles of the HJB and
		FP equations and shows how the general stochastic-to-deterministic route is applied in practice
		before tackling the more complex PV--battery case.		 
	\end{exmpl}

	\subsection{Numerical solution: Algorithm of multipliers}
	\noindent At this stage, we have a fully coupled optimality system. Its numerical solution therefore requires not only suitable time-space discretizations but also an iterative scheme. While \cite{BP} employs a method of multipliers based on the standard Lagrangian, here we propose an algorithm built on the Augmented Lagrangian formulation, leveraging the additional penalty term to improve conditioning and accelerate convergence.
	
	The algorithm consists of two coupled iterations: an outer loop that enforces the state constraints by updating the Lagrange multiplier $\boldsymbol{\mu}(s)$. At iteration $k$, 
	\begin{align}
		\nonumber
		\boldsymbol{\mu}^k(s) &= \rho (\partial \chi_{\mathcal{A}})_\lambda \left(\int_{\mathbb{R}^n} \mathbf{y}\varphi^k(s,\mathbf{y})\mathrm{d}\mathbf{y} + \lambda\boldsymbol{\mu}^{k-1}(s) \right) 
		\\	
		\label{eq:update_bfmu}
		&+ (1 - \rho) \boldsymbol{\mu}^{k-1}(s),
	\end{align}
	where $\rho \in \mathbb{R}$ is a relaxation parameter. Second, the inner loop solves an optimal control problem without state constraints. 
	
	A key innovation here, relative to \cite{BP}, is the insight that when $\boldsymbol{\mu}(s) = \mathbf{0}$, the inner loop collapses to a pure penalty scheme for the state constraints, with $\lambda$ playing the role of the penalty parameter. This observation motivates a two-stage strategy. First, a limited number of penalty-only iterations are performed by fixing $\boldsymbol{\mu}(s) = \mathbf{0}$, which guide the control towards a favorable descent direction while progressively enforcing the constraints. Subsequently, once a suitable regime has been reached, the method transitions to a standard Augmented Lagrangian framework, where both the multiplier $\boldsymbol{\mu}(s)$ and the penalty parameter are updated. This combined penalty-multiplier approach enhances both convergence speed and robustness with respect to the choice of the penalty parameter. A schematic summary of the resulting procedure is provided in Algorithm~\ref{alg:alg1}.
	
	\begin{algorithm}[!ht]
		\caption{Upgraded multipliers algorithm}\label{alg:alg1}
		\begin{algorithmic}
			\STATE Initialize control $\boldsymbol{\check{u}}^0(s, \mathbf{y})$.
			\STATE Set initial Lagrange multiplier $\boldsymbol{\mu}^0(s) = 0$. 
			\STATE Choose initial parameter $\lambda^0$ and target value $\lambda^*$.
			\STATE Choose convergence tolerances $\tau^*, \, \zeta^*$ and $\eta^*$.
			\STATE Set initial tolerance $\tau^0$.
			\STATE Set decreasing parameters $\omega_1, \, \omega_2 \in (0,1)$.
			\FOR{$k = 1, 2, 3, \dots$}
			\STATE Control $\boldsymbol{\check{u}}^{k-1}$ and Lagrange multiplier $\boldsymbol{\mu}^{k-1}$ are known.\STATE STEP 1: Compute $\varphi^k$ by solving FP equation, \eqref{eq:os11_1_laug}-\eqref{eq:os12_1_laug}, with $\boldsymbol{\check{u}}^{k-1}$.
			\STATE STEP 2: Compute $v^k$ and $\boldsymbol{\check{u}}^k$ by solving HJB equation, \eqref{eq:os12_1_laug}-\eqref{eq:os13_1_laug}, with $\boldsymbol{\check{u}}^{k-1}$, $\boldsymbol{\mu}^{k-1}$ and $\varphi^k$.
			\IF{$\|\boldsymbol{\check{u}}^k - \boldsymbol{\check{u}}^{k-1}\| \leq \tau^k$}
			\IF{$\|\boldsymbol{\check{u}}^k - \boldsymbol{\check{u}}^{k-1}\| \leq \tau^*$ \\
				\AND $\|\boldsymbol{\mu}^k - \boldsymbol{\mu}^{k-1}\| \leq \zeta^*$ \\
				\AND $\|\int_{\mathbb{R}^n} \mathbf{y}\varphi^k\mathrm{d}\mathbf{y} - \textrm{P}_\mathcal{A} \left(\int_{\mathbb{R}^n} \mathbf{y}\varphi^k\mathrm{d}\mathbf{y}\right)\| \leq \eta^*$ }
			\STATE \textbf{break}
			\ELSE
			\STATE Update Lagrange multiplier $\boldsymbol{\mu}^k$ as in \eqref{eq:update_bfmu}.
			\STATE Update parameter $\lambda$ and tighten tolerance $\tau$:
			\begin{align*}
				\tau^k &= \max\{\omega_1\tau^{k-1}, \tau^*\},
				\\
				\lambda^{k} &= \max\{\omega_2\lambda^{k-1}, \lambda^*\}.
			\end{align*}
			\ENDIF 
			\ELSE
			\STATE Take $\boldsymbol{\mu}^k = \boldsymbol{\mu}^{k-1}$, $\lambda^k = \lambda^{k-1}$ and $\tau^k = \tau^{k-1}$.
			\ENDIF
			\ENDFOR
		\end{algorithmic}
	\end{algorithm}

	\section{Spatial dimension reduction} \label{sect:IV_space_red}
	
	\noindent The development carried out so far requires solving both FP and HJB equations, each involving as many spatial variables as the dimension of the vector stochastic process defined by the original state equation \eqref{eq:SDE_system}. This leads to the well-known \textit{curse of dimensionality}, a term coined by Richard Bellman, which refers to the exponential increase in computational complexity as the number of state variables grows.
	
	Recently, novel approaches have been developed to overcome the curse of dimensionality. In particular, combining artificial neural networks with Monte Carlo methods to express semilinear PDE solutions as expectations of backward SDEs trajectories has shown great promise for solving HJB equations (see \cite{BE}). Although a full treatment of these cutting-edge techniques is outside the scope of this paper, in this section we introduce a simplified variant that still substantially mitigates the dimensionality challenge, while remaining fully consistent with the overall algorithmic framework presented.
	
	\paragraph{State decomposition and modeling assumptions}
	The central idea is to decompose the state process 
	\(\mathbf{X}^{t,\mathbf{x}}_s \in \mathbb{R}^n\) into two components \(\mathbf{X}^{t,\mathbf{x}_1}_{1,s} \in \mathbb{R}^{n_1}\) and \(\mathbf{X}^{t,\mathbf{x}_2}_{2,s} \in \mathbb{R}^{n_2}\), with \(n_1+n_2=n\):
	\[
	\mathbf{X}^{t,\mathbf{x}}_s=
	\begin{pmatrix}
		\mathbf{X}^{t,\mathbf{x}_1}_{1,s}\\
		\mathbf{X}^{t,\mathbf{x}_2}_{2,s}
	\end{pmatrix}.
	\]
	We identify the portion of the state directly affected by the control, \(\mathbf{X}^{t,\mathbf{x}_1}_{1,s}\), while the remaining part, \(\mathbf{X}^{t,\mathbf{x}_2}_{2,s}\), evolves independently of the control. We further require \(\mathbf{X}^{t,\mathbf{x}_1}_{1,s}\) and \(\mathbf{X}^{t,\mathbf{x}_2}_{2,s}\) to be statistically independent.
	
	Under suitable conditions, the uncontrolled block admits an additional factorization into independent lower-dimensional subprocesses:
	\[
	\mathbf{X}^{t,\mathbf{x}_2}_{2,s}=
	\begin{pmatrix}
		\mathbf{X}^{t,\mathbf{x}_{2,1}}_{2,1,s}\\
		\vdots\\
		\mathbf{X}^{t,\mathbf{x}_{2,q}}_{2,q,s}
	\end{pmatrix},
	\quad
	\mathbf{X}^{t,\mathbf{x}_{2,i}}_{2,i,s}\in\mathbb{R}^{n_{2,i}},\;\;
	\sum_{i=1}^{q} n_{2,i}=n_2.
	\]
	This setting enables a dimensionality reduction by replacing the original dynamics with several independent SDEs: one controlled SDE for \(\mathbf{X}^{t,\mathbf{x}_1}_{1,s}\) and \(q\) uncontrolled SDEs for the subprocesses \(\mathbf{X}^{t,\mathbf{x}_{2,i}}_{2,i,s}\), \(i=1,\dots,q\).

	Assume that the drift and diffusion fields decompose as
	\[
	\mathbf{b}(s,\mathbf{X}^{t,\mathbf{x}}_s,\boldsymbol{u}_s)=
	\begin{pmatrix}
		\mathbf{b}_1(s,\mathbf{X}^{t,\mathbf{x}_1}_{1,s},\boldsymbol{u}_s)\\
		\mathbf{b}_2(s,\mathbf{X}^{t,\mathbf{x}_2}_{2,s})
	\end{pmatrix},
	\quad
	\mathbf{\Sigma}(s,\mathbf{X}^{t,\mathbf{x}}_s)=
	\begin{pmatrix}
		\mathbf{\Sigma}_1(s,\mathbf{X}^{t,\mathbf{x}_1}_{1,s}) & \mathbf{0}\\
		\mathbf{0} & \mathbf{\Sigma}_2(s,\mathbf{X}^{t,\mathbf{x}_2}_{2,s})
	\end{pmatrix},
	\]
	and let \(\mathbf{B}_s=(\mathbf{B}_{1,s},\mathbf{B}_{2,s})^{\mathtt{T}}\). Suppose moreover that
	\[
	\mathbf{b}_2(s,\mathbf{X}^{t,\mathbf{x}_2}_{2,s})=
	\begin{pmatrix}
		\mathbf{b}_{2,1}(s,\mathbf{X}^{t,\mathbf{x}_{2,1}}_{2,1,s})\\
		\vdots\\
		\mathbf{b}_{2,q}(s,\mathbf{X}^{t,\mathbf{x}_{2,q}}_{2,q,s})
	\end{pmatrix},
	\quad
	\mathbf{\Sigma}_2(s,\mathbf{X}^{t,\mathbf{x}_2}_{2,s})
	=\operatorname{diag}\big(\mathbf{\Sigma}_{2,1},\dots,\mathbf{\Sigma}_{2,q}\big),
	\]
	with \(\mathbf{B}_{2,s}=(\mathbf{B}_{2,1,s},\dots,\mathbf{B}_{2,q,s})^{\mathtt{T}}\). Then the processes satisfy the uncoupled SDEs
	\begin{subequations}\label{eq:SDE_decomp_1}
		\begin{align}
			d\mathbf{X}^{t,\mathbf{x}_1}_{1,s} &=\mathbf{b}_1(s,\mathbf{X}^{t,\mathbf{x}_1}_{1,s},\boldsymbol{u}_s)\,ds +\mathbf{\Sigma}_1(s,\mathbf{X}^{t,\mathbf{x}_1}_{1,s})\,d\mathbf{B}_{1,s},
			\quad
			\mathbf{X}^{t,\mathbf{x}_1}_{1,t}=\mathbf{x}_1,\\
			d\mathbf{X}^{t,\mathbf{x}_{2,i}}_{2,i,s} &=\mathbf{b}_{2,i}(s,\mathbf{X}^{t,\mathbf{x}_{2,i}}_{2,i,s})\,ds +\mathbf{\Sigma}_{2,i}(s,\mathbf{X}^{t,\mathbf{x}_{2,i}}_{2,i,s})\,d\mathbf{B}_{2,i,s},
			\quad
			\mathbf{X}^{t,\mathbf{x}_{2,i}}_{2,i,t}=\mathbf{x}_{2,i},
			\;\; i=1,\dots,q.
		\end{align}
	\end{subequations}
	
	\paragraph{Objective and constraints}
	Although \(\mathbf{X}^{t,\mathbf{x}_1}_{1,s}\) and the subprocesses composing \(\mathbf{X}^{t,\mathbf{x}_2}_{2,s}\) evolve via uncoupled SDEs, they remain coupled at the \textit{optimization} level through the objective functional and the state constraints. The objective to be maximized is
	\begin{align}
		\nonumber
		J(t,\mathbf{x}_1,\mathbf{x}_{2,1},\dots,\mathbf{x}_{2,q},\boldsymbol{u})
		=\mathbb{E}\Bigg[&
		\int_t^T f\big(s,\mathbf{X}_{1,s}^{t,\mathbf{x}_1},
		\mathbf{X}_{2,1,s}^{t,\mathbf{x}_{2,1}},\dots,
		\mathbf{X}_{2,q,s}^{t,\mathbf{x}_{2,q}},\boldsymbol{u}_s\big)\,ds\\
		&\qquad\qquad\qquad\qquad
		+\,g\big(\mathbf{X}_{1,T}^{t,\mathbf{x}_1},
		\mathbf{X}_{2,1,T}^{t,\mathbf{x}_{2,1}},\dots,
		\mathbf{X}_{2,q,T}^{t,\mathbf{x}_{2,q}}\big)\Bigg].
		\label{eq:Fun_Obj_simpl}
	\end{align}
	A general formulation of state constraints is
	\[
	\mathbb{E}\big[\mathbf{G}(\mathbf{X}^{t,\mathbf{x}_1}_{1,s},
	\mathbf{X}_{2,1,s}^{t,\mathbf{x}_{2,1}},\dots,
	\mathbf{X}_{2,q,s}^{t,\mathbf{x}_{2,q}})\big]\in\mathcal{A},
	\, s\in[t,T],
	\]
	for a prescribed admissible set \(\mathcal{A}\subset\mathbb{R}^m\).
	
	\paragraph{Control structure and pdf factorization}
	Note that the dynamics in \eqref{eq:SDE_decomp_1} for the
	uncontrolled subprocesses do not involve \(\boldsymbol{u}_s\). Hence these components cannot be directly steered (e.g., by imposing bounds) at the path level. However, they may appear in the aggregated state constraints above, whose satisfaction depends indirectly on the control through \(\mathbf{X}^{t,\mathbf{x}_1}_{1,s}\).
	
	To recast the problem in a deterministic (FP-based) framework, we invoke the DPP and assume that an optimal feedback control exists depending only on the controlled state:
	\[
	\boldsymbol{u}_s=\boldsymbol{\check{u}}\big(s,\mathbf{X}^{t,\mathbf{x}_1}_{1,s}\big).
	\]
	Under this structure, \(\mathbf{X}^{t,\mathbf{x}_1}_{1,s}\) is independent of \(\mathbf{X}^{t,\mathbf{x}_2}_{2,s}\), and the subprocesses forming \(\mathbf{X}^{t,\mathbf{x}_2}_{2,s}\) are mutually independent. Therefore the joint density of \(\mathbf{X}^{t,\mathbf{x}}_s\) factorizes as
	\[
	\varphi(s,\mathbf{y}_1,\mathbf{y}_{2,1},\dots,\mathbf{y}_{2,q})
	=\varphi_1(s,\mathbf{y}_1)\,\varphi_2(s,\mathbf{y}_2)
	=\varphi_1(s,\mathbf{y}_1)\,\prod_{i=1}^q \varphi_{2,i}(s,\mathbf{y}_{2,i}),
	\]
	with \(\mathbf{y}_1\in\mathbb{R}^{n_1}\), \(\mathbf{y}_2\in\mathbb{R}^{n_2}\),
	\(\mathbf{y}_{2,i}\in\mathbb{R}^{n_{2,i}}\). Here, \(\varphi_1\) and \(\varphi_2\) are the pdfs of \(\mathbf{X}_{1,s}^{t,\mathbf{x}_1}\) and \(\mathbf{X}_{2,s}^{t,\mathbf{x}_2}\), and \(\varphi_{2,i}\) is the pdf of \(\mathbf{X}_{2,i,s}^{t,\mathbf{x}_{2,i}}\).
	
	\paragraph{Implications for FP equations}
	The above factorization yields a corresponding decomposition of the FP equation into independent components: one governing \(\varphi_1\), that depends on the feedback control, and \(q\) equations for \(\varphi_{2,1},\dots,\varphi_{2,q}\), that are control-free. Since the coefficients \(\mathbf{b}_{2,i}\) and \(\mathbf{\Sigma}_{2,i}\) are given for \(i=1,\dots,q\), the FP equations for the uncontrolled subprocesses can be solved independently \textit{a priori}, and then used in the reduced-dimension optimal control problem.
	
	\paragraph{On the scope of the decomposition}
	Before proceeding with the deterministic reformulation and the derivation of the optimality system, we highlight a limitation of the methodology. A reduction of the state dimension is \textit{not} feasible if the optimal feedback control depends on \(\mathbf{X}^{t,\mathbf{x}_1}_{1,s}\) and on \textit{every} component of \(\mathbf{X}^{t,\mathbf{x}_2}_{2,s}\). In that case, the SDE for \(\mathbf{X}^{t,\mathbf{x}_1}_{1,s}\) would involve \(\mathbf{X}^{t,\mathbf{x}_2}_{2,s}\), the marginal density \(\varphi_1\) would depend on both blocks, and the factorization of \(\varphi\) would no longer hold. One must then revert to the general (non-reduced) framework.
	
	A meaningful reduction is still possible when the optimal feedback depends on \(\mathbf{X}^{t,\mathbf{x}_1}_{1,s}\) and on a strict subset of the uncontrolled components, say \(\mathbf{X}^{t,\mathbf{x}_{2,i}}_{2,i,s}\), \(i=1,\dots,p\) with \(p<q\). In that situation, we regroup the state as
	\[
	\mathbf{X}^{t,\mathbf{x}}_s=
	\begin{pmatrix}
		\bar{\mathbf{X}}^{t,\mathbf{x}_1}_{1,s}\\[2pt]
		\bar{\mathbf{X}}^{t,\mathbf{x}_2}_{2,s}
	\end{pmatrix},
	\quad
	\bar{\mathbf{X}}^{t,\mathbf{x}_1}_{1,s}=
	\begin{pmatrix}
		\mathbf{X}^{t,\mathbf{x}_{1}}_{1,s}\\
		\mathbf{X}^{t,\mathbf{x}_{2,1}}_{2,1,s}\\
		\vdots\\
		\mathbf{X}^{t,\mathbf{x}_{2,p}}_{2,p,s}
	\end{pmatrix},
	\quad
	\bar{\mathbf{X}}^{t,\mathbf{x}_2}_{2,s}=
	\begin{pmatrix}
		\mathbf{X}^{t,\mathbf{x}_{2,p+1}}_{2,p+1,s}\\
		\vdots\\
		\mathbf{X}^{t,\mathbf{x}_{2,q}}_{2,q,s}
	\end{pmatrix}.
	\]
	Only the FP equations associated with the remaining uncontrolled blocks \(\varphi_{2,i}\), \(i=p+1,\dots,q\), can be solved in advance as they are control-free). The others enter the controlled block and must be handled within the reduced optimality system.
	
	\paragraph{Illustration via PV--battery problem}
	Consider the PV--battery example of Section~\ref{sect:II_I_PV_problem} with state and control
	\[
	\mathbf{X}^{t,\mathbf{x}}_s=
	\begin{pmatrix}
		Z_s\\ \Pi_s\\ \mathcal{E}_s
	\end{pmatrix},
	\quad
	u_s=P_{\mathrm{bat},s},
	\]
	where the control acts only on the SDE for \(\mathcal{E}_s\). Under the structural assumptions on \(\mathbf{b}\) and \(\mathbf{\Sigma}\), a first decomposition is
	\[
	\mathbf{X}^{t,\mathbf{x}_1}_{1,s}=\mathcal{E}_s,
	\quad
	\mathbf{X}^{t,\mathbf{x}_2}_{2,s}=
	\begin{pmatrix}
		Z_s\\ \Pi_s
	\end{pmatrix},
	\]
	with a further split \(\mathbf{X}^{t,\mathbf{x}_{2,1}}_{2,1,s}=Z_s\), \(\mathbf{X}^{t,\mathbf{x}_{2,2}}_{2,2,s}=\Pi_s\). Then the pdfs of \(Z_s\) and \(\Pi_s\) can be precomputed and the optimal feedback depends only on the battery energy, \(P_{\mathrm{bat},s}=\check{P}_{\mathrm{bat}}(s,\mathcal{E}_s)\). If the control authority is too limited, one may enlarge the controlled block to include prices,
	\[
	\mathbf{X}^{t,\mathbf{x}_1}_{1,s}=
	\begin{pmatrix}
		\Pi_s\\ \mathcal{E}_s
	\end{pmatrix},
	\quad
	\mathbf{X}^{t,\mathbf{x}_2}_{2,s}=Z_s,
	\]
	so that \(P_{\mathrm{bat},s}=\check{P}_{\mathrm{bat}}(s,\Pi_s,\mathcal{E}_s)\) is richer, while only the pdf of \(Z_s\) remains available in advance.
	
	\paragraph{Design trade-off}
	The choice of decomposition must balance computational tractability against control expressiveness. Enlarging the controlled block improves responsiveness by incorporating richer state information, but reduces the part that can be precomputed. Conversely, a smaller controlled block accelerates computation and simplifies implementation, at the possible cost of reduced performance. Importantly, this dimension reduction strategy does not introduce numerical approximation errors; instead, it yields a simplified, non-equivalent optimal control problem defined on a reduced state space. Consequently, any performance differences should be interpreted as arising from the limited state information available to the controller.
		
	A systematic quantitative assessment of how state reduction affects control accuracy is challenging in this setting, as it depends on problem-specific dynamics, objectives, and constraints, and is therefore beyond the scope of this work. Nevertheless, the numerical results include a comparative evaluation based on the objective functional, illustrating how the chosen decomposition influences control performance in practice when applied to PV-battery systems.

	\paragraph{Deterministic reformulation and notation}
	We now reformulate the stochastic control problem as a deterministic (FP-based) control problem. To simplify notation, for the second block of the state we write functions in terms of a single variable \(\mathbf{y}_2\in\mathbb{R}^{n_2}\) rather than the components \(\mathbf{y}_{2,i}\in\mathbb{R}^{n_{2,i}}\), \(i=1,\dots,q\). We take \(\varphi_1(s,\mathbf{y}_1)\) as the state variable. Accordingly, the state equation is its associated FP equation, and the control is the feedback function \(\boldsymbol{\check{u}}(s,\mathbf{y}_1)\). The objective functional becomes
	\begin{align*}
		J(\boldsymbol{\check{u}}) &= \int_t^T \!\Bigg(\int_{\mathbb{R}^{n_1}} 
		\bigg( \int_{\mathbb{R}^{n_2}} f\big(s,\mathbf{y}_1,\mathbf{y}_2,\boldsymbol{\check{u}}(s,\mathbf{y}_1)\big)\,
		\varphi_2(s,\mathbf{y}_2)\,d\mathbf{y}_2\bigg)\,\varphi_1(s,\mathbf{y}_1)\,d\mathbf{y}_1\Bigg)\,ds \\
		&\quad + \int_{\mathbb{R}^{n_1}} \left(\int_{\mathbb{R}^{n_2}} g(\mathbf{y}_1,\mathbf{y}_2)\,\varphi_2(T,\mathbf{y}_2)\,d\mathbf{y}_2\right)\,
		\varphi_1(T,\mathbf{y}_1)\,d\mathbf{y}_1,
	\end{align*}
	where the expectation is evaluated using \(\varphi_1\) and the precomputed \(\varphi_2\). In particular, when a finer factorization is available one may set
	\(\varphi_2=\prod_{i=1}^q \varphi_{2,i}\). Analogously, state constraints read
	\[
	\int_{\mathbb{R}^{n_1}}\!\left( \int_{\mathbb{R}^{n_2}} \mathbf{G}(\mathbf{y}_1,\mathbf{y}_2)\,\varphi_2(s,\mathbf{y}_2)\,d\mathbf{y}_2\right)\,
	\varphi_1(s,\mathbf{y}_1)\,d\mathbf{y}_1 \in \mathcal{A},\, s\in[t,T].
	\]
	
	\paragraph{Optimality system}
	The first-order conditions follow the approach of \cite{BP}, with the difference that state constraints are enforced via an Augmented Lagrangian. The resulting system is
	{\allowdisplaybreaks
		\begin{subequations} \label{eq:opt_system_decomp}
			\begin{align}
				\nonumber
				&\frac{\partial \varphi_1}{\partial s} (s,\mathbf{y}_1)+\mathrm{div}\Big(\mathbf{b}_1(s,\mathbf{y}_1, \boldsymbol{\check{u}}(s,\mathbf{y}_1))\varphi_1(s,\mathbf{y}_1)\Big) 
				\\*
				\label{eq:os11_1_laug_Y}
				& \hspace*{1.45cm} -\frac{1}{2}\mathrm{div} \mathbf{div}\Big(\mathbf{\Sigma}_1 (s,\mathbf{y}_1)\mathbf{\Sigma}_1(s,\mathbf{y}_1)^\mathtt{T}\varphi_1(s,\mathbf{y}_1)\Big)=0,
				\\*
				\label{eq:os12_1_laug_Y}
				&\varphi_1(t,\mathbf{y}_1) = \varphi_{1,\mathbf{x}_1}(\mathbf{y}_1),\ \mathbf{y}_1\in\mathbb{R}^{n_1},
				\\
				\nonumber
				&\frac{\partial v}{\partial s}(s,\mathbf{y}_1) +\max_{\boldsymbol{u}\in\mathrm{\mathbf{U}}}\bigg\{\mathbf{b}_1(s,\mathbf{y}_1,\boldsymbol{u})\cdot\mathrm{grad} \, v(s,\mathbf{y}_1)
				\\*
				\nonumber
				&\hspace*{2.4cm} +\int_{\mathbb{R}^{n_2}} f(s,\mathbf{y}_1, \mathbf{y}_2, \boldsymbol{u}) \varphi_2(s,\mathbf{y}_2)\mathrm{d}\mathbf{y}_2\bigg\}
				\\*
				\nonumber
				&\hspace*{1.4cm}+\frac{1}{2}\mathbf{\Sigma}_1 (s,\mathbf{y}_1)\mathbf{\Sigma}_1 (s,\mathbf{y}_1)^\mathtt{T}\cdot \mathrm{Hess}\,(v)(s,\mathbf{y}_1) 
				\\*
				\nonumber
				&\hspace*{1.4cm}- (\partial \chi_{\mathcal{A}})_\lambda \Bigg(\int_{\mathbb{R}^{n_1}} \bigg( \int_{\mathbb{R}^{n_2}} \mathbf{G}(\mathbf{y}_1, \mathbf{y}_2) \varphi_2(s,\mathbf{y}_2) \mathrm{d}\mathbf{y}_2 \bigg)
				\\*
				\nonumber
				&\hspace*{3.8cm}\varphi_1 (s,\mathbf{y}_1)\mathrm{d}\mathbf{y}_1 + \lambda \boldsymbol{\mu}(s)\Bigg)
				\\*
				\label{eq:os13_1_laug_Y}
				&\hspace*{1.6cm} \cdot \int_{\mathbb{R}^{n_2}} \mathbf{G}(\mathbf{y}_1, \mathbf{y}_2) \varphi_2(s,\mathbf{y}_2) \mathrm{d}\mathbf{y}_2 = 0,
				\\*
				\label{eq:os14_1_laug_Y}
				&  v(T,\mathbf{y}_1)=\int_{\mathbb{R}^{n_2}}g(\mathbf{y}_1, \mathbf{y}_2)\varphi_2(T,\mathbf{y}_2)\mathrm{d}\mathbf{y}_2,\ \mathbf{y}_1\in\mathbb{R}^{n_1},
				\\
				\nonumber
				&\boldsymbol{\mu}(s) = (\partial \chi_{\mathcal{A}})_\lambda \Bigg(\int_{\mathbb{R}^{n_1}}\left( \int_{\mathbb{R}^{n_2}} \mathbf{G}(\mathbf{y}_1, \mathbf{y}_2) \varphi_2(s,\mathbf{y}_2) \mathrm{d}\mathbf{y}_2 \right) 
				\\*
				\label{eq:os15_1_laug_Y}
				&\hspace{2.9cm}\varphi_1 (s,\mathbf{y}_1)\mathrm{d}\mathbf{y}_1 + \lambda \boldsymbol{\mu}(s)\Bigg),
				\\
				\nonumber
				& \boldsymbol{\check{u}}(s,\mathbf{y}_1)=\underset{\boldsymbol{u}\in\mathrm{\mathbf{U}}}{\text{arg max}}\bigg\{\mathbf{b}_1(s,\mathbf{y}_1,\boldsymbol{u})\cdot\mathrm{grad}\, v(s,\mathbf{y}_1)
				\\*
				\label{eq:os16_1_laug_Y}
				&\hspace*{2.7cm}+ \int_{\mathbb{R}^{n_2}} f(s,\mathbf{y}_1, \mathbf{y}_2, \boldsymbol{u}) \varphi_2(s,\mathbf{y}_2)\mathrm{d}\mathbf{y}_2\bigg\}.
			\end{align}
		\end{subequations}
		
		\paragraph{Computational remarks}
		In practice, one first computes \(\varphi_2\) (or \(\{\varphi_{2,i}\}_{i=1}^q\)) offline. Then the reduced optimality system in dimension \(n_1\) is solved (e.g., by Algorithm~\ref{alg:alg1}). This strategy is typically far more efficient than tackling the original system in dimension \(n>n_1\), while it preserves the essential coupling to the uncontrolled block through the precomputed \(\varphi_2\).
		
		From a practical standpoint, if the state of the problem can be decomposed into two independent stochastic processes, the proposed methodology not only accelerates the numerical resolution but also reduces the number of state variables directly influencing the control. }

	\section{Numerical results} \label{sect:V_results}
	\noindent In this section, we return to the stochastic control problem for the grid-connected PV system with batteries and apply the concepts developed in the preceding sections. Additionally, two alternatives for reducing the space dimension are proposed and numerical results are presented to illustrate the effectiveness of the approach.

	\subsection{PV power plant}
	\noindent First, let us identify the components of the PV plant problem within the proposed general formulation:
	
	\begin{itemize}
		\item The state is $\mathbf{X}_s = (Z_s, \Pi_s, \mathcal{E}_s)^\mathtt{T}$ and the optimality system has space dimension $n=3$.
		\item The control is $u_s = P_{\text{bat},s}$.
		\item Vector field $\mathbf{b}$ is
		\begin{equation*}
			\mathbf{b} (s, \mathbf{y},u)=
			\begin{pmatrix}
				\kappa_{Z}\left(\theta_{Z}(s)
				-y_1\right)
				\\
				\\
				\kappa_{\Pi}\left(\theta_{\Pi}(s)+\frac{1}{\kappa_{\Pi}}\frac{d\theta_\Pi}{ds}(s)-y_2\right)
				\\
				\\
				-u
			\end{pmatrix}.
		\end{equation*}
		\item Tensor field $\mathbf{\Sigma}$ is 
		\begin{equation*}
			\mathbf{\Sigma} (s, \mathbf{y})=
			\begin{pmatrix}
				\sigma_{Z,Z} y_1 (1 - y_1) & 0 & 0
				\\
				\\
				\sigma_{\Pi,Z} y_2 & \sigma_{\Pi,\Pi} y_2 &0
				\\
				\\
				0 & 0 & \sigma_{\mathcal{E},\mathcal{E}} y_3
			\end{pmatrix}.
		\end{equation*}
		\item The Brownian motion is $\mathbf{B}_s = (B_{Z,s}, B_{\Pi,s}, B_{\mathcal{E},s})^\mathtt{T}$.
		\item $f(s, \mathbf{y}, u) = y_2 \left(\hat{P}(I_{CS}(s)y_1) + u\right). $
		\item $g(\mathbf{y}) = c_T y_3$, where $c_T$ (\euro/MWh) is a given parameter.
		\item $\mathrm{U} = \{u : P_{min}\leq u \leq P_{max}\}.$
		\item $\mathcal{A} = \{ \mathbf{y} \in \mathbb{R}^3 : \mathcal{E}_{min}\leq y_3 \leq \mathcal{E}_{max}\}.$ 
	\end{itemize}
	
	It should be noted that this problem presents a special characteristic: an analytical expression for the feedback optimal control, $\check{P}_{\text{bat}}(s, \mathbf{y})$, can be derived. Indeed, by denoting $\mathbf{p}(s, \mathbf{y}) = \mathrm{grad}\, v(s,\mathbf{y})$, it is verified that
	\begin{align*}
		\check{P}_{\text{bat}}(s, \mathbf{y})&=
		\underset{u\in\mathrm{U}}{\text{arg max}} \left[\mathbf{b}(s,\mathbf{y},u) \cdot \mathbf{p}(s, \mathbf{y}) + f(s, \mathbf{y}, u) \right] 
		\\
		&=\underset{u\in\mathrm{U}}{\text{arg max}} \bigg[b_1(s,\mathbf{y},u) p_1(s, \mathbf{y}) + b_2(s,\mathbf{y},u) p_2(s, \mathbf{y})
		\\
		& \hspace*{1.7cm} - u p_3(s, \mathbf{y})+ y_2\left(\hat{P}(I_{CS}(s)y_1) + u\right) \bigg]
		\\
		&=\underset{u\in\mathrm{U}}{\text{arg max}} \left[(-p_3(s, \mathbf{y}) + y_2)u \right].
	\end{align*}
	If $\mbox{meas}\left\lbrace(s,\mathbf{y})\in [t,T]\times\mathbb{R}^3: \, -p_3(s, \mathbf{y})+y_2 =0\right\rbrace=0,$ then the optimal control is bang-bang and given by
	\begin{align}
		\label{eq:feedback_opt_ctrl}
		\check{P}_{\text{bat}}(s, \mathbf{y}) &= P_{\text{bat},s} = \begin{cases}
			P_{min} & \mbox{if } -p_3+y_2 <0,
			\\
			P_{max} & \mbox{if } -p_3 +y_2>0.
		\end{cases}
	\end{align}
	As a consequence, the HJB equation \eqref{eq:os13_1_laug}, \eqref{eq:os14_1_laug} can be written as the following semilinear PDE:
	\begin{subequations} \label{eq:HJB_PV_system}
		\begin{align}
			\nonumber
			\frac{\partial v}{\partial s}(s,\mathbf{y}) &+ \kappa_{Z}\left(\theta_{Z}(s)-y_1\right)\frac{\partial v}{\partial y_1}(s,\mathbf{y})
			\\
			\nonumber
			&+ \kappa_{\Pi}\left(\theta_{\Pi}(s)+
			\frac{1}{\kappa_{\Pi}}\frac{d\theta_\Pi}{ds}(s)-y_2\right)
			\frac{\partial v}{\partial y_2}(s,\mathbf{y})
			\\
			\nonumber
			& +\frac{1}{2}\sigma_{Z,Z}^2 y_1^2 (1-y_1)^2 \frac{\partial^2 v}{\partial y_1^2}(s,\mathbf{y})
			\\
			\nonumber
			&+\sigma_{Z,Z}y_1(1-y_1)\sigma_{\Pi,Z} y_2 \frac{\partial^2 v}{\partial y_1\partial y_2}(s,\mathbf{y})
			\\
			\nonumber
			&+\frac{1}{2}(\sigma_{\Pi,Z} ^2 +\sigma_{\Pi,\Pi}^2)y_2^2\frac{\partial^2v}{\partial y_2^2}(s,\mathbf{y}) +\frac{1}{2}\sigma_{\mathcal{E},\mathcal{E}}^2y_3^2\frac{\partial^2v}{\partial y_3^2}(s,\mathbf{y}) 
			\\
			\nonumber
			&+\Psi\left(y_2,\frac{\partial v}{\partial y_3}(s,\mathbf{y})\right) + y_2 \hat{P}(I_{CS}(s)y_1)
			\\
			\label{eq:HJB_PV}
			& -(\partial \chi_{\mathcal{A}})_\lambda \bigg(\int_{\mathbb{R}^n}\mathbf{y}\varphi (s,\mathbf{y})\,\mathrm{d}\mathbf{y} + \lambda \boldsymbol{\mu}(s)\bigg) \cdot \mathbf{y}=0,
			\\
			\label{eq:HJB_end_PV}
			v(T,\mathbf{y})&=c_T y_3,
		\end{align}
	\end{subequations}
	where $\Psi$ is the non-linear function defined as
	\begin{equation*}
		\Psi(z,q) =\begin{cases}
			P_{min}(-q+z) & \mbox{if}\ -q+z<0,
			\\
			P_{max}(-q+z) & \mbox{if}\ -q+z>0.
		\end{cases}
	\end{equation*}
	
	\subsection{Dimension-reduction techniques}
	\noindent Following the lines of Section~\ref{sect:IV_space_red}, the decomposition of the state vector into two components in the form 
	\begin{align*}
		\mathbf{X}_{1,s}^{t,\mathbf{x}_1} = \mathcal{E}_s, 
		\quad 
		\mathbf{X}_{2,s}^{t,\mathbf{x}_2} = \left(\begin{array}{c}
			Z_s \\ 
			\Pi_s
		\end{array}\right),
	\end{align*} 
	leads to a reduction of the spatial dimension to $n_1=1$. Consequently, the optimization problem can be simplified significantly, especially since optimal control and value function now depend only on~$\mathcal{E}_s$.
	
	The corresponding decompositions of vector and tensor fields, $\mathbf{b}$ and $\mathbf{\Sigma}$, can be directly deduced to get an optimality system with space dimension $n_1 = 1$. Additionally, the optimization problem characterizing the feedback control can be solved analytically. Indeed, denoting $p(s, y_{1,1}): = \mathrm{grad}\, v(s, y_{1,1})$, the optimal control in feedback form is \footnote{In general, $y_{i,j}$ denotes the $j$-th component of vector $\mathbf{y}_i$}
	\begin{align}
		\label{eq:feedback_opt_ctrl_decomp_1}
		\check{P}_{\text{bat}}(s,y_{1,1}) &= \begin{cases}
			P_{min} & \mbox{if } -p(s, y_{1,1})+ \mathrm{E}[y_{2,2}] <0,
			\\
			P_{max} & \mbox{if } -p(s, y_{1,1}) + \mathrm{E}[y_{2,2}]>0.
		\end{cases}
	\end{align}
	
	It is important to highlight that this simplification restricts decision-making to the battery's state of charge. In practical applications, such a limited perspective may not be adequate for participation in day-ahead electricity market auctions or for real-time control. To address this limitation, we propose an alternative decomposition that incorporates electricity price information into the control strategy of the PV system. Specifically, 
	\begin{align*}
		\mathbf{X}_{1,s}^{t,\mathbf{x}_1} = \left( \begin{array}{c}
			\Pi_s \\ 
			\mathcal{E}_s 
		\end{array}\right), 
		\quad
		\mathbf{X}_{2,s}^{t,\mathbf{x}_2} = Z_s.
	\end{align*}
	
	Now, the optimality system to be solved has spatial dimension $n_1 =2$, and then optimal control and value function depend on both $\Pi_s$ and $\mathcal{E}_s$. Again, an explicit expression for the feedback control can be derived. It corresponds to
	\begin{align}
		\label{eq:feedback_opt_ctrl_decomp_2}
		\check{P}_{\text{bat}}(s,\mathbf{y}_{1}) = \begin{cases}
			P_{min} & \mbox{if } -p_2(s,\mathbf{y}_{1}) + y_{1,1} <0,
			\\
			P_{max} & \mbox{if } -p_2(s,\mathbf{y}_{1}) + y_{1,1}>0,
		\end{cases}
	\end{align}
	where $\mathbf{p}(s,\mathbf{y}_{1}) = \mathrm{grad}\, v(s, \mathbf{y}_{1})$ and hence $p_2 = \partial v/\partial \mathcal{E}$. In this case, HJB equation \eqref{eq:os13_1_laug}, \eqref{eq:os14_1_laug} reads
	\begin{subequations}\label{eq:HJB_PV_system_decomp_2}
		\begin{align}
			\nonumber
			&\frac{\partial v}{\partial s}(s, \mathbf{y}_1) + \kappa_{\Pi}\left(\theta_{\Pi}(s)+\frac{1}{\kappa_{\Pi}}\frac{d\theta_\Pi}{ds}(s)-y_{1,1}\right) \frac{\partial v}{\partial y_{1,1}}(s, \mathbf{y}_1)
			\\
			\nonumber
			&+ \frac{1}{2}\sigma_{\Pi,\Pi}^2y_{1,1}^2\frac{\partial^2v}{\partial y_{1,1}^2}(s,\mathbf{y}_1) + \frac{1}{2}\sigma_{\mathcal{E},\mathcal{E}}^2y_{1,2}^2\frac{\partial^2v}{\partial y_{1,2}^2}(s,\mathbf{y}_1)
			\\
			\nonumber
			&+\Psi\Big(y_{1,1},\frac{\partial v}{\partial y_{1,2}}(s,\mathbf{y}_1)\Big) 
			+ \int_{\mathbb{R}} y_{1,1}\hat{P}(I_{CS}(s)y_{2,1}) \varphi_2(s, y_{2,1})\,\mathrm{d} y_{2,1} 
			\\
			\label{eq:HJB_PV_decomp_2}
			&- (\partial \chi_{\mathcal{A}})_\lambda \bigg(\int_{\mathbb{R}^{n_1}}\mathbf{y}_1\varphi_1 (s,\mathbf{y}_1)\, \mathrm{d} \mathbf{y}_1 + \lambda \boldsymbol{\mu}(s)\bigg) \cdot \mathbf{y}_1 = 0,
			\\
			\label{eq:HJB_end_PV_decomp_2}
			&v(T,\mathbf{y}_1)=c_T y_{1,2},
		\end{align}
	\end{subequations}
	where 
	\begin{equation*}
		\Psi(z,q) :=\left\{\begin{array}{cc}
			P_{min}(-q+z) & \mbox{if }\ - q+z<0,
			\\
			P_{max}(-q+z) & \mbox{if }\ -q+z>0.
		\end{array}\right.
	\end{equation*}
	
	It can be shown that the solution to this equation, $v(s, \mathbf{y}_1)$, also satisfies the expectation of equation \eqref{eq:HJB_PV_system} with respect to the stochastic process CSI, $Z_s$. Firstly, note that since $v(s, \mathbf{y}_1)$ does not depend on $Z_s$ (i.e., on $y_{2,1}$), its derivatives with respect to $y_{2,1}$ vanish. Secondly, by multiplying equation \eqref{eq:HJB_PV_system} by $\varphi_2$ and then integrating over $\mathbb{R}$ with respect to $Z_s$ (i.e., to $y_{2,1}$), we can compute the expectation and derive the equation for $v(s, \mathbf{y}_1)$, which results in \eqref{eq:HJB_PV_system_decomp_2}.

	\subsection{Comparative benchmarking methods}
	
	\noindent In order to provide a comparative benchmarking of the proposed methodology, an MPC method and two rule-based methods are considered. These benchmarks serve to contextualize the performance of the proposed optimization-based methodology in terms of both economic performance and computational cost.

	\subsubsection{Rule-based methods}
	\noindent
	Rule-based approaches, although simpler than optimization-based methods, are widely used in the literature as baseline tools for PV--battery operation due to their lower computational requirements. An important distinction from the proposed DPP framework is that rule-based methods operate in discrete time. Therefore, a time step $\Delta s$ must be explicitly introduced.
	
	The first strategy is based on defining two price thresholds, $\Pi_{min}$ and $\Pi_{max}$. At each time step, the control action, namely, the battery charging or 
	discharging power $P_{\text{bat},s}$, is selected according to the following rule (recall the sign convention, $P_{\text{bat},s}$ is positive when the battery is discharging and negative when it is charging): 
	\begin{equation*}
		P_{\text{bat},s} =
		\begin{cases}
			P_{min}, & \mbox{if } \Pi_s < \Pi_{min} \mbox{ and } \mathcal{E}_s - \Delta s P_{min} < \mathcal{E}_{max}, \\
			P_{max}, & \mbox{if } \Pi_s > \Pi_{max} \mbox{ and } \mathcal{E}_s - \Delta s P_{max} > \mathcal{E}_{min}, \\
			0, & \mbox{otherwise}.
		\end{cases}
	\end{equation*}
	This rule ensures charging during low-price periods and 
	discharging during high-price periods, subject to battery limits. A third case, $P_{\text{bat},s} = 0$, corresponds to idle operation and prevents rapid switching between charging and discharging when prices fluctuate around the threshold values.

	The second strategy relies on a time-of-use (TOU) tariff division of the day. Two sets of time intervals, \textit{off-peak} and \textit{peak} hours, are defined in advance based on the historical trend of the electricity price process. The rule prescribes charging during off-peak hours, discharging during peak hours, and idle operation otherwise:
	\begin{equation*}
		P_{\text{bat},s} =
		\begin{cases}
			P_{min}, & \mbox{if } s \in \mbox{off-peak hours and } \mathcal{E}_s - \Delta s P_{min} < \mathcal{E}_{max}, \\
			P_{max}, & \mbox{if } s \in \mbox{peak hours and } \mathcal{E}_s - \Delta s P_{max} > \mathcal{E}_{min}, \\
			0, & \text{otherwise}.
		\end{cases}
	\end{equation*}
	
	To assess these rules within the stochastic framework, the stochastic state dynamics are simulated using the Euler--Maruyama scheme and at each time step the corresponding rule-based method is applied. The expected revenues are then evaluated through a Monte Carlo procedure over multiple realizations of state and control.

	\subsubsection{MPC method}
	
	In addition to the proposed continuous-time Dynamic Programming approach, a stochastic MPC strategy is considered as an optimization-based benchmark. MPC is widely used in energy management applications due to its ability to incorporate system constraints and updated information over a time horizon.
	
	At initial step, MPC solves the optimal control problem over the prediction
	horizon T, but applies only the control action along a small interval $[t, t + \Delta s]$, then repeats the process for time interval $[t + \Delta s, T]$ taking the real state of the system at time $t$ as initial condition for the model.
	Let $s_k$ denote the current decision time. At each stage $k$, the controller solves a finite-horizon stochastic optimal control problem by maximizing the objective function 
	\begin{equation*}
		J_k(\boldsymbol{u}) = 
		\mathbb{E}\left[
		\int_{s_k}^{T} f(s, \mathbf{X}_{s}^{s_k, \mathbf{x}}, \boldsymbol{u}_s)\,ds
		+ g(\mathbf{X}_{T}^{s_k, \mathbf{x}})
		\right],
	\end{equation*}
	subject to the PV plant stochastic system dynamics and the control constraints previously introduced. As in the proposed framework, state constraints are imposed on the expected value of the state variables, i.e.,
	\[
	\mathbb{E}[\mathcal{E}_{s}] \in [\mathcal{E}_{min}, \mathcal{E}_{max}], \qquad s \in [s_k, T].
	\]
	
	In order to solve the stochastic optimal control problem at each MPC stage, stochastic Pontryagin Maximum Principle (SPMP) is applied and state constraints are handled using an Augmented Lagrangian approach. Therefore, uncertainty is accounted for through Monte Carlo sampling of the stochastic dynamics. Finally, the MPC stage $k$ yields a control sequence over the horizon $[s_k, T]$ that maximizes the objective function while ensuring the state constraints imposed over the state expectation.
	
	Since the MPC formulation optimizes a single state trajectory rather than a feedback control law, the state used to initialize the next MPC stage is taken as the expected value of the state resulting from the application of the first control action. This choice allows the MPC scheme to follow an average trajectory comparable to the expected-state evolution obtained under the Dynamic Programming approach.
	
	The resulting MPC algorithm can be summarized as follows:
	\begin{itemize}
		\item At time $s_k$, solve the finite-horizon stochastic optimal control problem with state constraints enforced in expectation.
		\item Apply the first control action of the optimal sequence.
		\item Update the state using its expected value and repeat the procedure at $s_{k+1}$.
	\end{itemize}

	\subsection{Data}
\begin{figure}[!t]
		\centering
		\includegraphics[width=3.in]{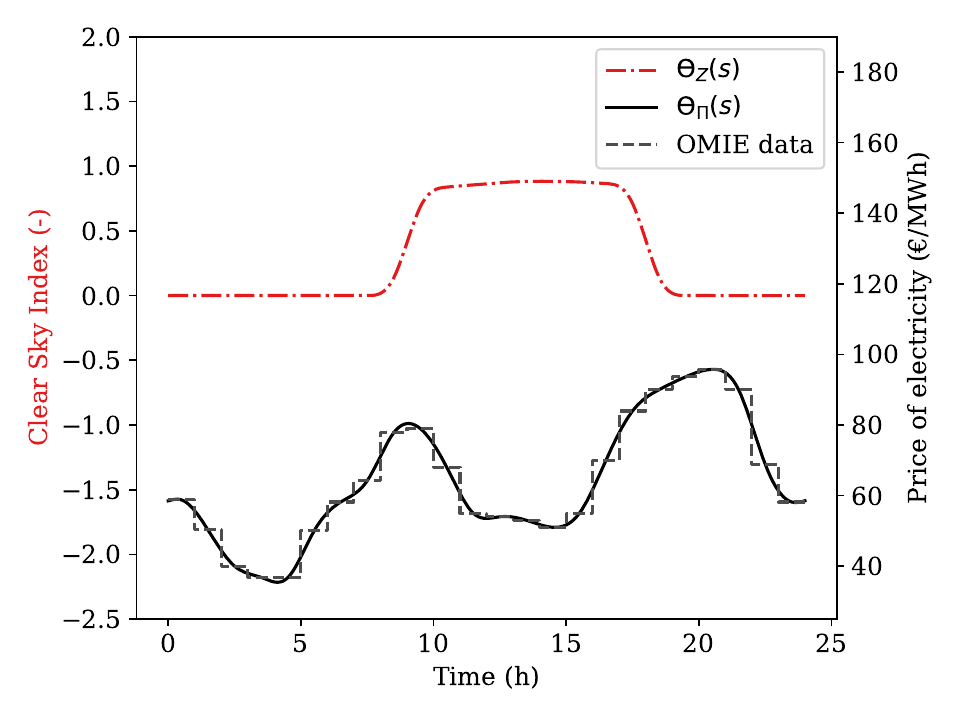}
		\caption{Functions $\theta_{Z}(s)$ and $\theta_{\Pi}(s)$.}
		\label{fig:data_theta}
	\end{figure}

	\noindent Data used for the numerical resolution of the PV plant problem are listed below. Some of them are taken from \cite{BGGK, JHD}:
	\begin{itemize}
		\item Time interval to be considered is $s\in(0, 24)$~h.
		\item Mean reverting parameters (in h$^{-1}$): $\kappa_{Z} = 0.75 $ and $\kappa_{\Pi}= 0.04$.
		\item Functions $I_{CS}(s)$ and $\theta_Z(s)$ have been derived from historical data for the first month of 2022, with a frequency of 1 hour. The Clear Sky Irradiance of every day of the month, $I_{CS} (s)$, has been computed by evaluating the expression~\eqref{eq:I_CS}, and then daily average was estimated based on these values. To obtain $\theta_Z (s)$, first the average daily GHI, $I(s)$, is estimated in a similar manner. Plots of the mean irradiances $I_{CS}(s)$ and $I(s)$ are shown in Fig.~\ref{fig:exp_irrad}. Then, $\theta_{Z}(s)$ is computed using definition \eqref{eq:Z_s}.
		\item Function $\theta_{\Pi}(s)$ has been taken from the Iberian market (MIBEL) as the hourly electricity price for December 21, 2022 (see \cite{OMIE}). Specifically, we have taken
		a cubic spline interpolation of the step function defined by the hourly price data. Fig.~\ref{fig:data_theta} depicts curves $\theta_{Z}(s)$ and $\theta_{\Pi}(s)$.
		\item Constant parameters in $\mathbf{\Sigma}$ (in h$^{-1/2}$): $\sigma_{Z,Z} = 0.2$, $\sigma_{\Pi,\Pi} =0.075$, $\sigma_{\mathcal{E},\mathcal{E}} = 0.1$ and $\sigma_{\Pi,Z} = 0.$ This parameter is set to zero to allow for the proper decomposition of the matrix $\mathbf{\Sigma}$ in the two dimension reduction techniques discussed earlier. In other scenarios, it may take a positive value to account for the influence of the renewable energy factor on electricity prices, as proposed in \cite{BP}.
		\item The initial conditions of the SDE model, $\mathbf{x}=(Z^0, \Pi^0, \mathcal{E}^0)^\mathtt{T}$, are taken as a normal random vector with mean vector and matrix variance given by
		\begin{align*}
			\bar{\mathbf{X}} = \left(\begin{array}{c}
				\theta_Z(0) \\
				\theta_\Pi (0) \\
				2
			\end{array}\right) \text{ and } 
			V_0 = \begin{pmatrix}
				0.1 & 0 & 0 \\
				0 & 8.560 & 0 \\
				0 & 0 & 0.01
			\end{pmatrix},
		\end{align*}
		respectively. The joint pdf of this random vector is taken as initial condition in FP equation, namely $\varphi_\mathbf{x}$.
		\item Parameters for PV power, $\hat{P}$, in \eqref{eq:PV_power}: $A_{PV} = 7500$ m$^2$, $\eta = 0.8$ and $\beta = 0.2382 \pi$~rad. Solar inclination angle $\gamma(s)$ has been taken as the daily average of historical data of the first month of 2022.
		\item Bounds for the power delivered by the battery: $P_{min}=-1$~MW and $P_{max} = 1$ MW.
		\item Bounds for the state constraint, namely, for the energy content of the battery: $\mathcal{E}_{min} = 0$ MWh and $\mathcal{E}_{max}=4$~MWh.
		\item Terminal cost: $c_T=\theta_\Pi(T)$~\euro/MWh. 
	\end{itemize}
	
	Additional data for rule-based methods include:
	\begin{itemize}
		\item Price thresholds: $\Pi_{min} = 65$~\euro/MWh and $\Pi_{max} = 75$~\euro/MWh.
		\item Off-peak hours and peak hours are defined according to Spanish electricity tariff schedule. Specifically, for one day off-peak hours extend from $00{:}00$ to $08{:}00$, while peak hours go from $10{:}00$ to $14{:}00$ and from $18{:}00$ to $22{:}00$.
		\item Number of realizations for Monte Carlo approach: $10000$.
	\end{itemize}
	
	\subsection{Results}
	\noindent To numerically solve the PV plant optimization problem, both the FP and HJB equations are discretized in time using the implicit Euler scheme, and in space using first-order finite element methods. Depending on the spatial dimension of the problem, interval, triangular, or tetrahedral meshes are employed. Additionally, the computational procedure follows Algorithm~\ref{alg:alg1}.
	
	The chosen computational domain for the state variables is $(-5, 5) \times (-40, 220) \times (-30, 40)$, ensuring a sufficiently large region to capture the solution's support so that the solution of the FP equations preserves the property of having integral equal 1. Results presented here correspond to a structured mesh consisting of $ \mathbf{N} = (21, 66, 71)^\mathtt{T}$ nodes for each spatial dimension. However, as discussed later, various mesh resolutions have been tested to assess convergence and accuracy. For the temporal discretization, time horizon is set to the interval $(0,24)$ hours, representing one day, with a uniform time step of $\Delta s = 0.5$~h. 
	
\begin{figure}[!t]
		\centering
		\includegraphics[width=3.in]{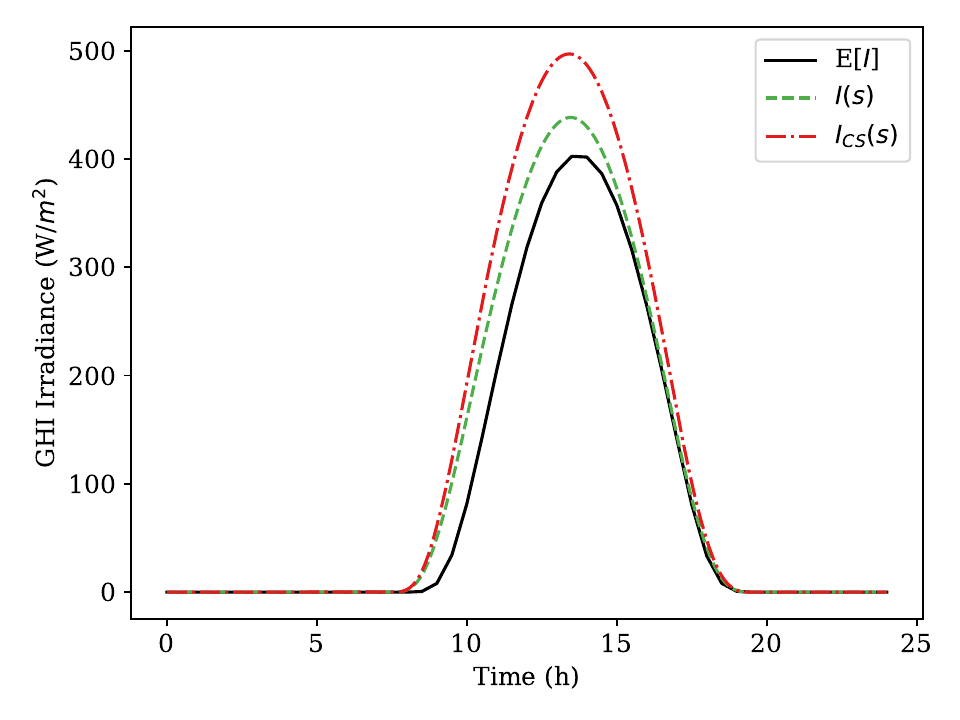}
		\caption{Expectation of $I_s$ and functions $I_{CS}(s)$ and $I(s)$.}
		\label{fig:exp_irrad}
	\end{figure}

\begin{figure}[!t]
		\centering
		\includegraphics[width=3.in]{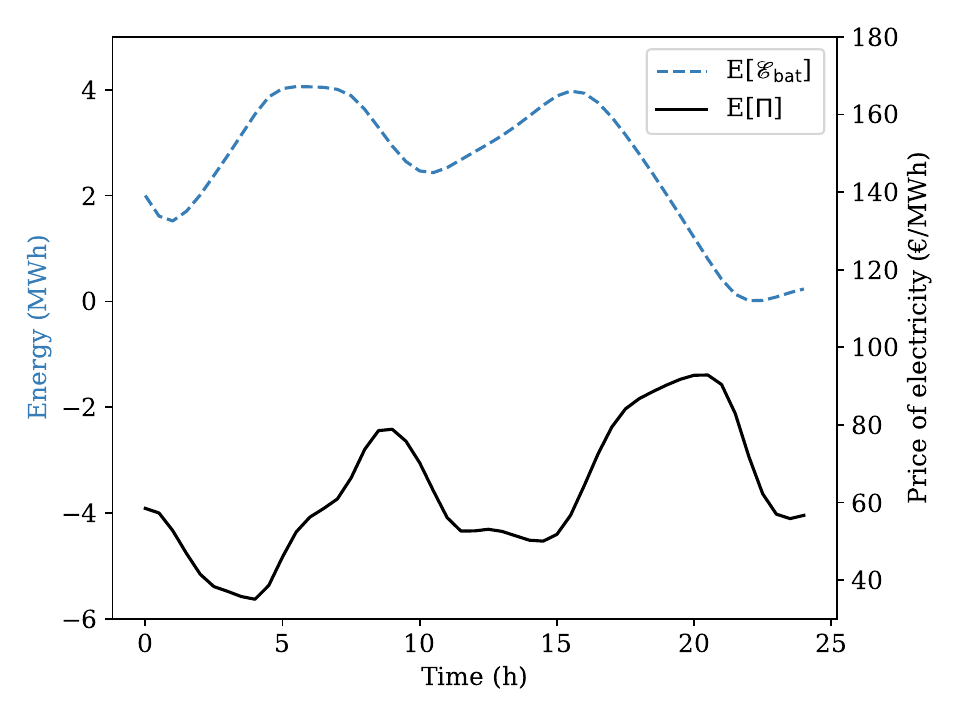}
		\caption{Expectation of $\Pi_s$ and $\mathcal{E}_s$.}
		\label{fig:exp_EPi}
	\end{figure}

	We begin by analyzing the optimal results of the full three-dimensional formulation of the PV plant problem. Fig.~\ref{fig:exp_irrad} displays the expected value of the solar irradiance $I_s$, along with the reverting profile $I(s)$ and the Clear Sky Irradiance $I_{CS}(s)$. The expectation of $I_s$ is computed based on the expected value of the first state variable using expression~\eqref{eq:Z_s}. As shown in the figure, irradiance is naturally concentrated during daylight hours and remains below the clear-sky benchmark. Similarly, Fig.~\ref{fig:exp_EPi} shows the time evolution of the expected values of the remaining state variables: electricity price and battery's energy level. As observed, the expected battery energy, $\mathrm{E}[\mathcal{E}_s]$, remains within the prescribed bounds of $0$ and $4$ MWh, indicating that the state constraints are correctly enforced by the algorithm. Moreover, the battery’s charging and discharging patterns align clearly with low- and high-price periods, respectively, underscoring the economic logic of the optimal control strategy.
	
	
\begin{figure}[!t]
		\centering
		\includegraphics[width=3.in]{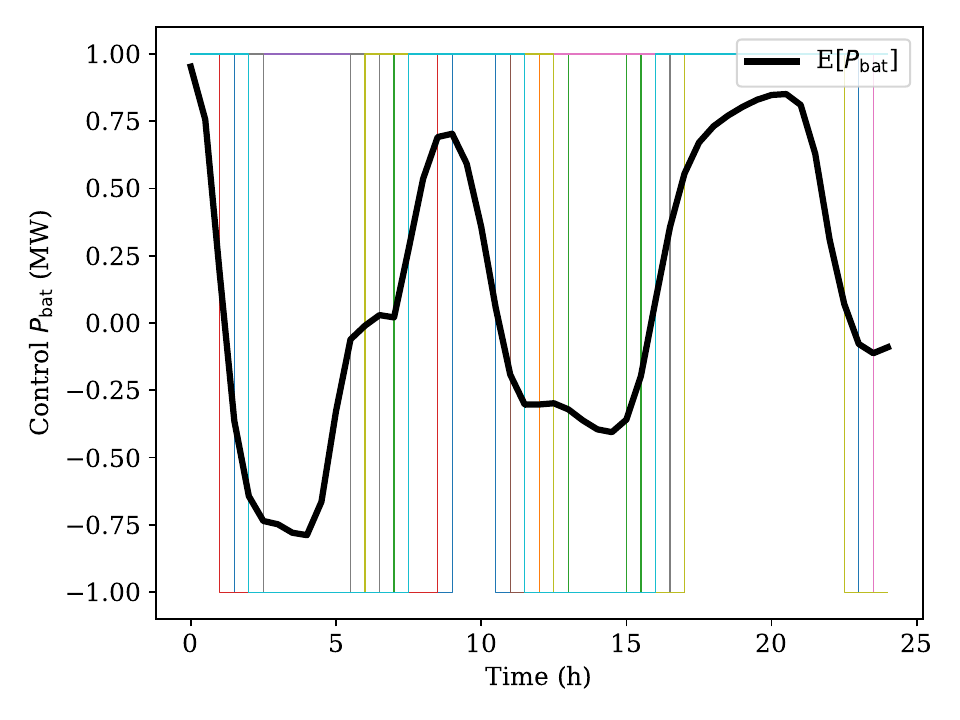}
		\caption{Expectation of $P_{\text{bat},s}$ and 10 realizations.}
		\label{fig:exp_Pbat}
	\end{figure}

	To illustrate the control behavior, Fig.~\ref{fig:exp_Pbat} presents not only the expected value of the battery power, $\mathrm{E}[P_{\text{bat},s}]$, but also ten sample trajectories of this stochastic process. These realizations clearly exhibit the characteristic bang-bang behavior of the optimal control. Moreover, the control constraints are consistently satisfied across all trajectories.
	
\begin{figure}[!t]
		\centering
		\includegraphics[width=3.in]{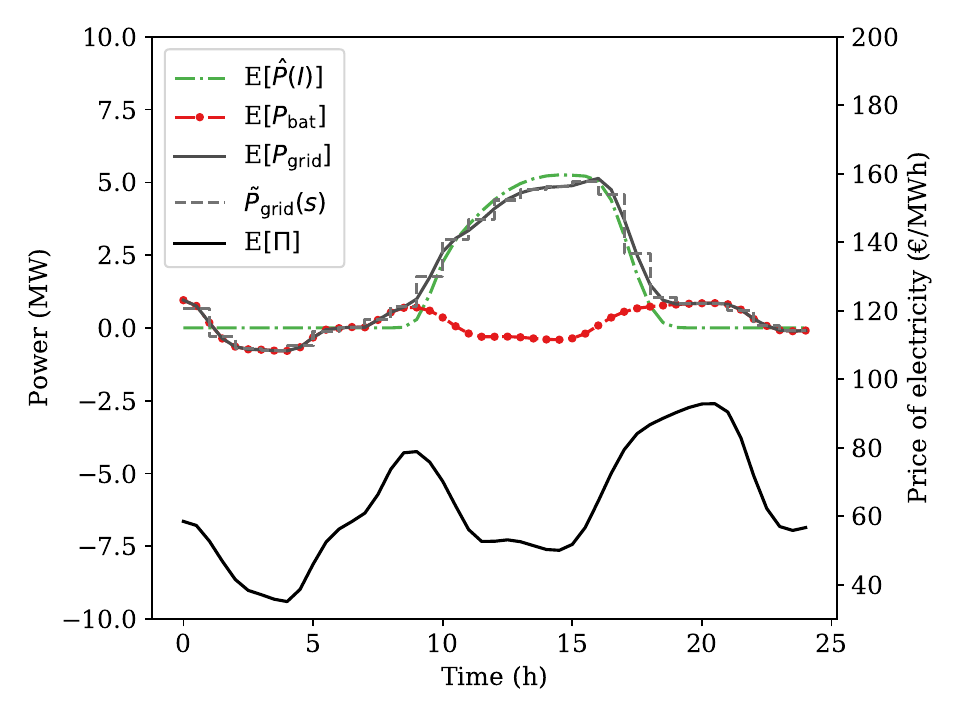}
		\caption{Expectation of $P_{\text{bat},s}$, $P_{\text{solar},s}$, $P_{\text{grid},s}$ and $\Pi_s$. Optimal power bid at day-ahead electricity market auctions is included (\textit{capacity firming)}.}
		\label{fig:exp_powers_Pi}
	\end{figure}

	Fig.~\ref{fig:exp_powers_Pi} compares the expected values of the power generated by the PV system, $\mathrm{E}[P_{\text{solar},s}]$, the power supplied by the battery, $\mathrm{E}[P_{\text{bat},s}]$, and the net power exchanged with the grid, $\mathrm{E}[P_{\text{grid},s}]$, alongside the expected electricity price $\mathrm{E}[\Pi_s]$. This representation highlights how generation, storage operation, and grid interaction are jointly coordinated in response to evolving price signals. In particular, the battery’s charging and discharging actions complement PV generation to take advantage of favorable market conditions. Additionally, the bids that the PV producer submits to the day-ahead electricity market auctions are also displayed, $\tilde{P}_{\text{grid}}(s)$. These bids are typically required to be piecewise constant over hourly intervals, a practice known as \textit{capacity firming}. To generate these bids, the approach described in equation \eqref{eq:exp_u} is applied and then averaged on an hourly basis. From an operational perspective, the results reveal a clear economic pattern: during periods of low solar generation and low electricity prices, the plant tends to purchase energy from the grid to charge the battery. Conversely, when prices rise, the producer prioritizes exporting PV generation and may additionally discharge the battery, thereby increasing the net power injected into the grid and maximizing revenue.
	
	
\begin{figure}[!t]
		\centering
		\includegraphics[width=3.in]{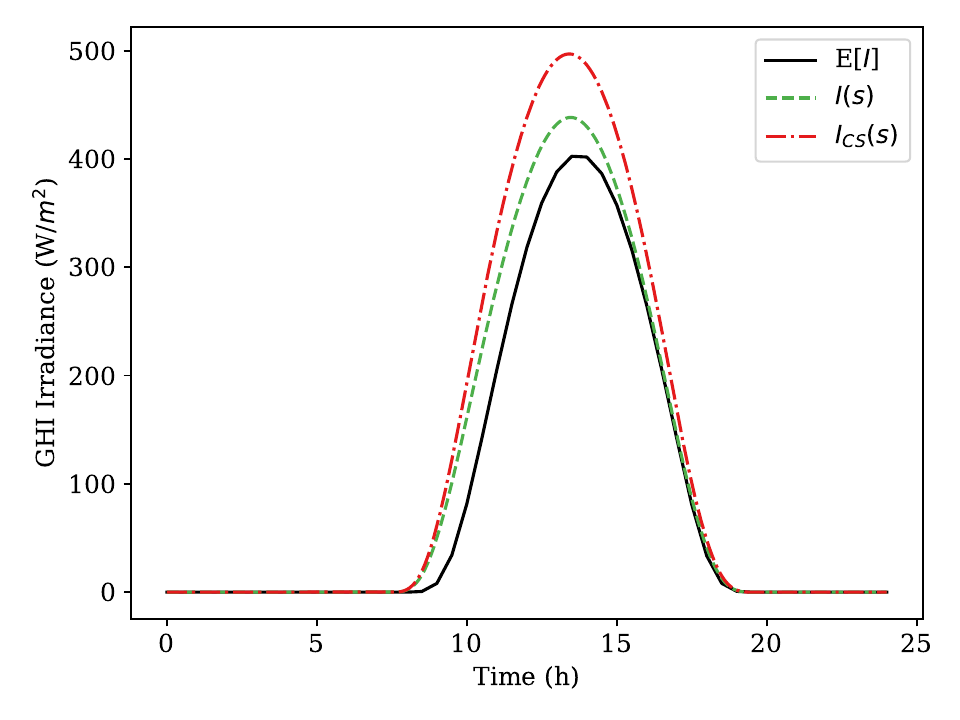}
		\caption{Expectation of $I_s$ and functions $I_{CS}(s)$ and $I(s)$ for 1D case.}
		\label{fig:exp_irrad_1D}
	\end{figure}

\begin{figure}[!t]
		\centering
		\includegraphics[width=3.in]{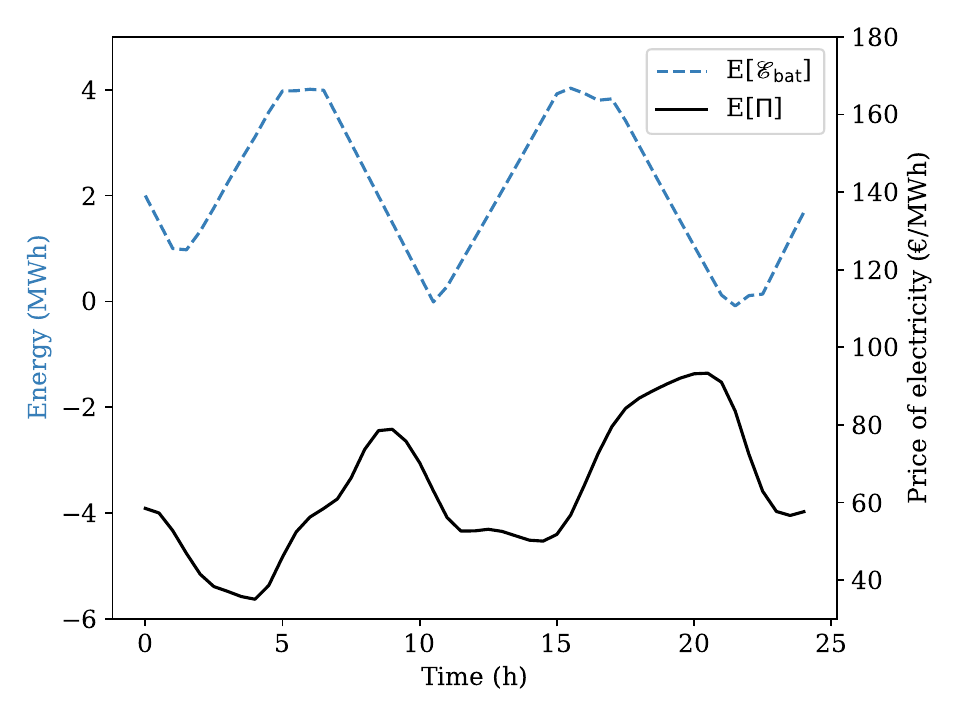}
		\caption{Expectation of $\Pi_s$ and $\mathcal{E}_s$ for 1D case.}
		\label{fig:exp_EPi_1D}
	\end{figure}

\begin{figure}[!t]
		\centering
		\includegraphics[width=3.in]{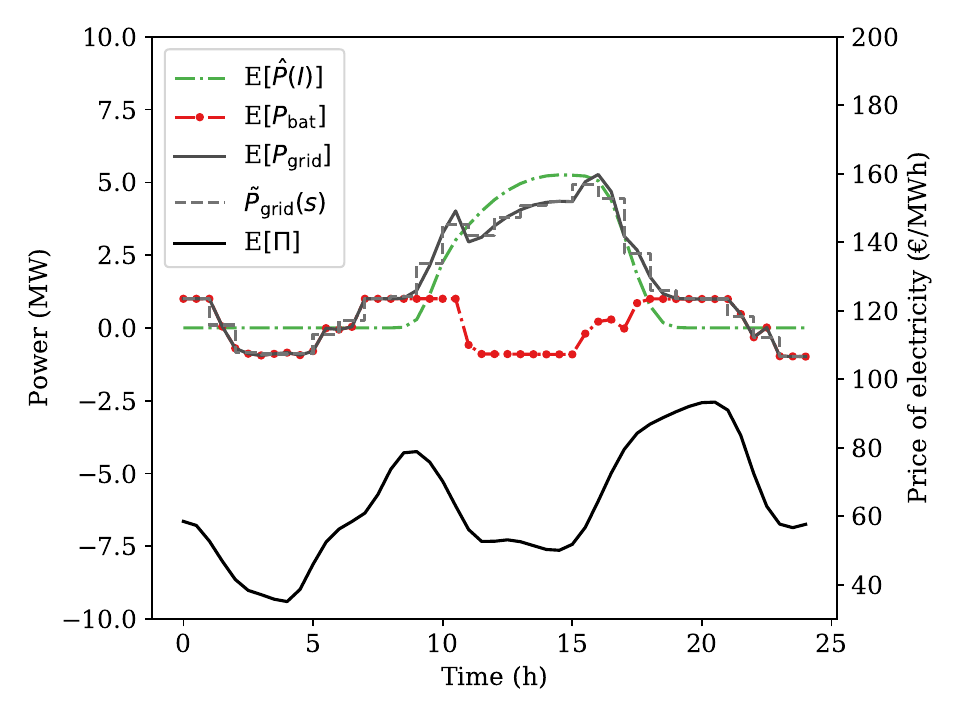}
		\caption{Expectation of $P_{\text{bat},s}$, $P_{\text{solar},s}$, $P_{\text{grid},s}$ and $\Pi_s$ together with optimal power bid $\tilde{P}_{\text{grid}}(s)$ for 1D case.}
		\label{fig:exp_powers_Pi_1D}
	\end{figure}

	Next, we present the numerical results obtained by solving the PV plant problem using the decomposition which reduces the system to one spatial dimension (1D). Fig.~\ref{fig:exp_irrad_1D} shows the expected solar irradiance, $\mathrm{E}[I_s]$, and functions $I_{CS}(s)$ and $I(s)$. Fig.~\ref{fig:exp_EPi_1D} depicts the expected evolution of the electricity price and the battery's energy content. Finally, Fig.~\ref{fig:exp_powers_Pi_1D} presents the time evolution of the expected control, $\mathrm{E}[P_{\text{bat},s}]$, the solar power $\mathrm{E}[P_{\text{solar},s}]$, the power exchanged with the grid $\mathrm{E}[P_{\text{grid},s}]$, and the electricity price $\mathrm{E}[\Pi_s]$. It also includes the hourly constant bids at the electricity market, $\tilde{P}_{\text{grid}}(s)$. As expected, the behavior of solar irradiance and electricity prices remains essentially unchanged relative to the three-dimensional case, since these state variables are not directly influenced by the control strategy. In contrast, noticeable differences arise in the evolution of the battery energy content and the charge-discharge power. Although all state and control constraints remain satisfied, the battery energy more frequently approaches its upper and lower bounds, indicating tighter constraint saturation. This behavior reflects the reduced expressiveness of the control policy in the 1D setting, where decisions rely on a smaller subset of state variables. Nevertheless, the qualitative operating pattern remains consistent with the full 3D formulation: the battery tends to charge during low-price periods and discharge when prices are high, confirming that the reduced model preserves the main economic logic of the optimal strategy, although with a less flexible control response.
	
\begin{figure}[!t]
		\centering
		\includegraphics[width=3.in]{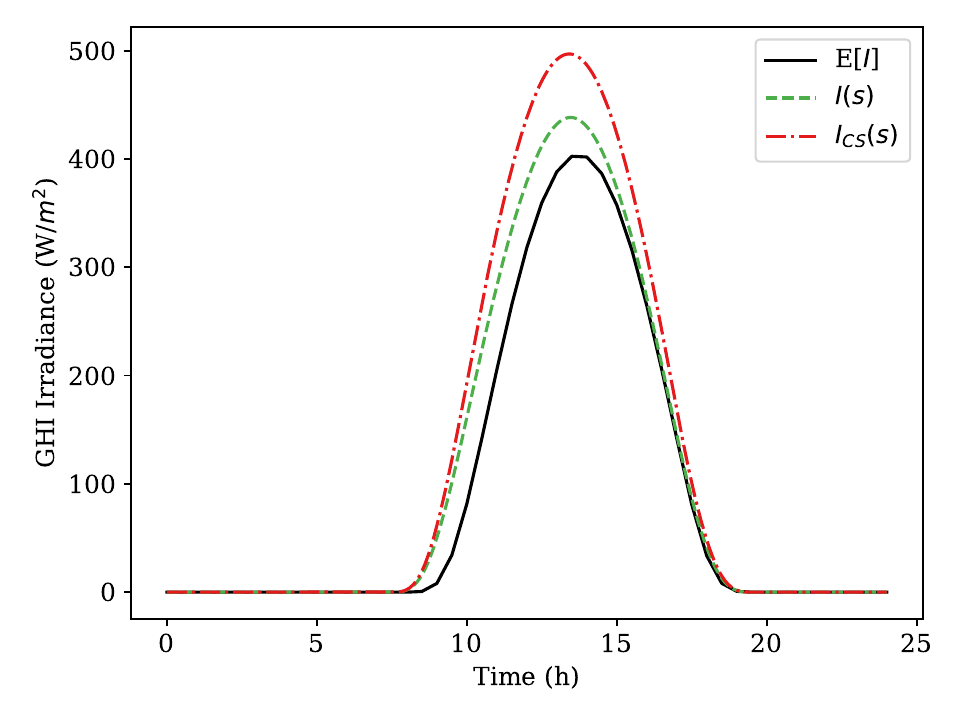}
		\caption{Expectation of $I_s$ and functions $I_{CS}(s)$ and $I(s)$ for 2D case.}
		\label{fig:exp_irrad_2D}
	\end{figure}

\begin{figure}[!t]
		\centering
		\includegraphics[width=3.in]{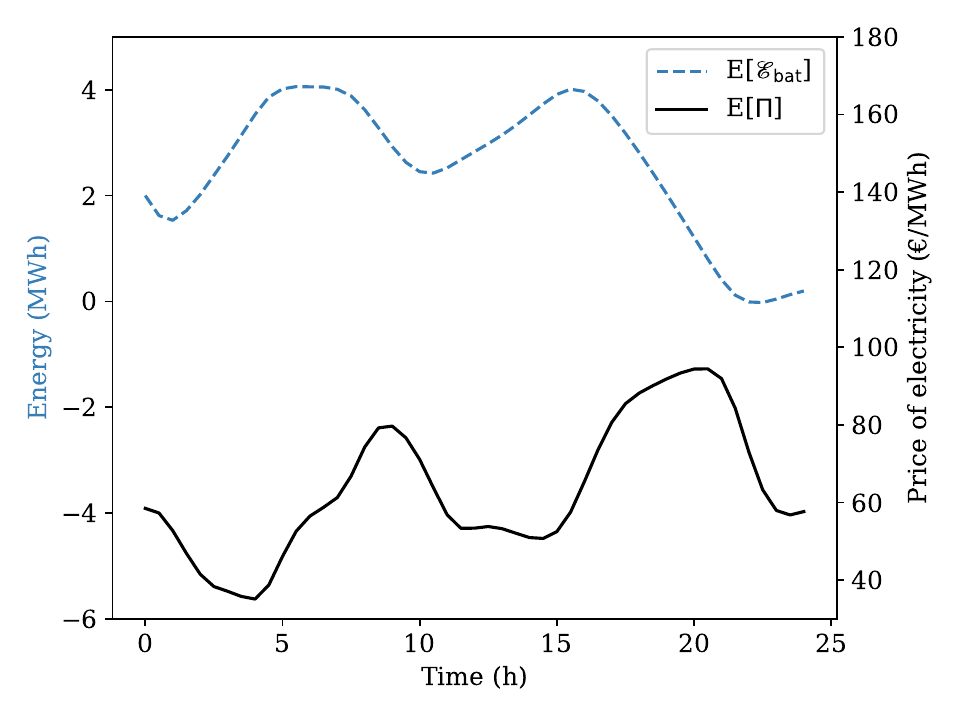}
		\caption{Expectation of $\Pi_s$ and $\mathcal{E}_s$ for 2D case.}
		\label{fig:exp_EPi_2D}
	\end{figure}

\begin{figure}[!t]
		\centering
		\includegraphics[width=3.in]{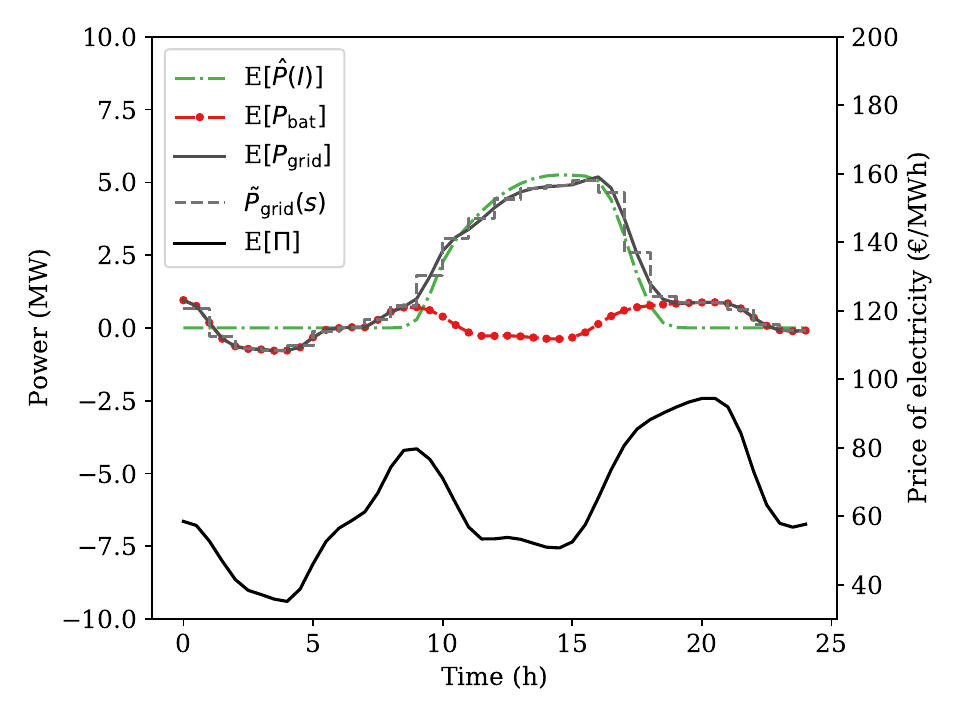}
		\caption{Expectation of $P_{\text{bat},s}$, $P_{\text{solar},s}$, $P_{\text{grid},s}$ and $\Pi_s$, and optimal power bid $\tilde{P}_{\text{grid}}(s)$ for 2D case.}
		\label{fig:exp_powers_Pi_2D}
	\end{figure}

	For the decomposition leading to a 2D problem, we present analogous results to those of the 1D case. These are shown in Fig.~\ref{fig:exp_irrad_2D}, which displays the expected solar irradiance $\mathrm{E}[I_s]$; in Fig.~\ref{fig:exp_EPi_2D}, showing the expectations of the electricity price $\mathrm{E}[\Pi_s]$ and the battery energy content $\mathrm{E}[\mathcal{E}_s]$; and in Fig.~\ref{fig:exp_powers_Pi_2D}, illustrating the expected values of the various power flows along with $\mathrm{E}[\Pi_s]$. In contrast to the 1D case, incorporating electricity prices into the feedback control law significantly improves the coordination between battery charge-discharge actions and price signals. As a result, the battery energy trajectory exhibits less frequent saturation of its bounds, and the overall behavior more closely resembles that of the full 3D formulation while still satisfying all constraints. This demonstrates the benefit of enriching the control policy with relevant information, yet maintaining a substantial reduction in computational complexity, since only 2D PDEs must be solved.

\begin{figure}[!t]
		\centering
		\includegraphics[width=4.in]{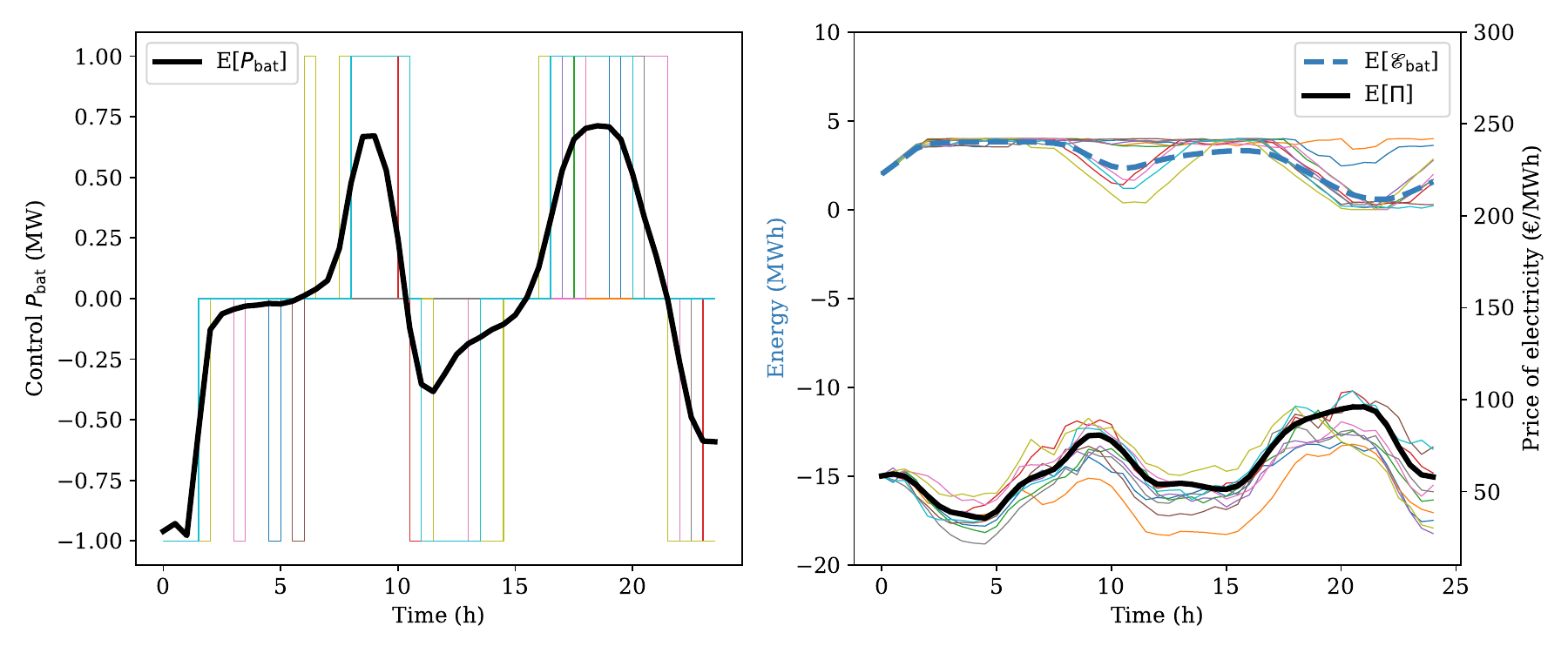}
		\caption{Expectation of $\mathcal{E}_s$, $\Pi_s$ and $P_{\text{bat},s}$ along with 10 realizations via price threshold method.}
		\label{fig:exp_price_threshold}
	\end{figure}

\begin{figure}[!t]
		\centering
		\includegraphics[width=4.in]{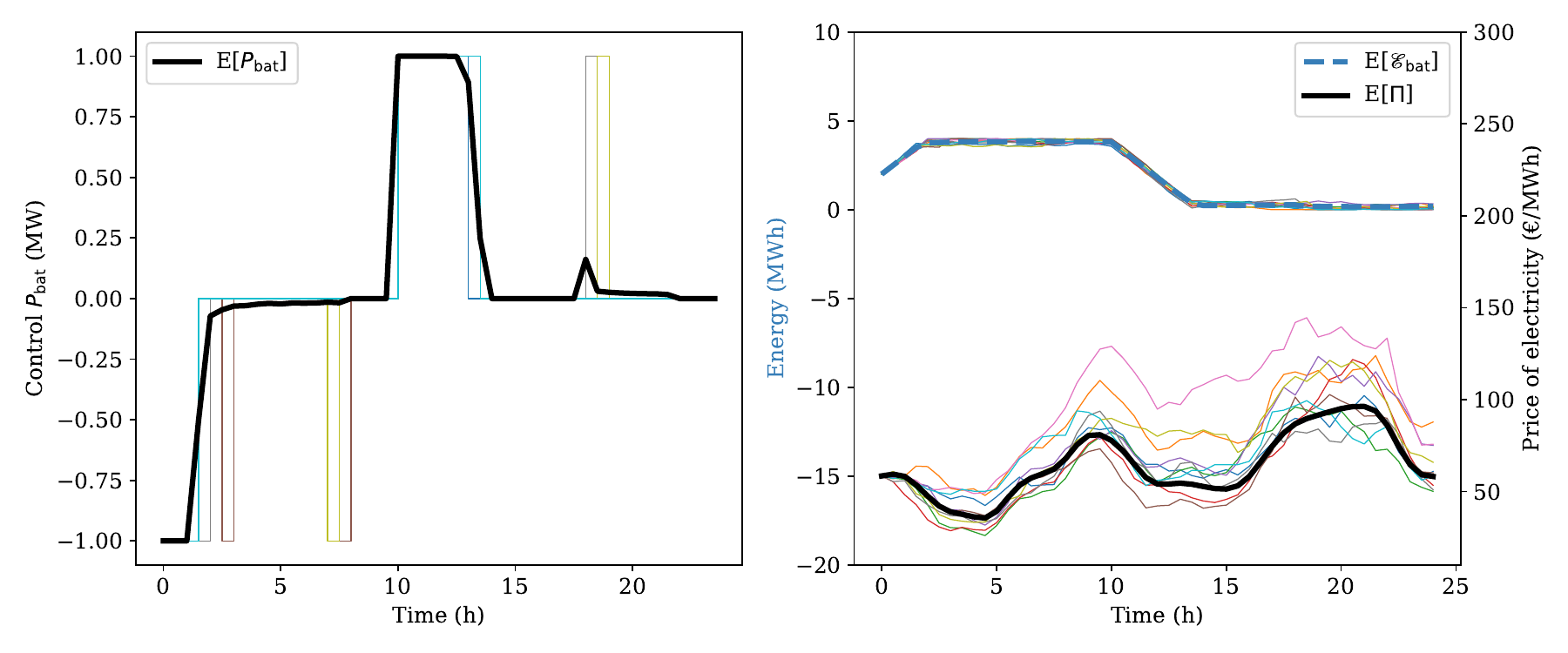}
		\caption{Expectation of $\mathcal{E}_s$, $\Pi_s$ and $P_{\text{bat},s}$ along with 10 realizations via TOU method.}
		\label{fig:exp_TOU}
	\end{figure}

	Regarding the rule-based methods, Fig.~\ref{fig:exp_price_threshold} illustrates the expected battery power, $\mathrm{E}[P_{\text{bat},s}]$, energy content, $\mathrm{E}[\mathcal{E}_s]$ and electricity price, $\mathrm{E}[\Pi_s]$ together with 10 sample realizations of each variable obtained under the price-threshold rule. Analogously, Fig.~\ref{fig:exp_TOU} presents the corresponding results for the TOU strategy. In both cases, state constraints are enforced directly during the simulation and are therefore satisfied. The price-threshold rule exhibits a qualitative behavior that is broadly consistent with the proposed methodology, in the sense that charging and discharging actions tend to align with low- and high-price periods, respectively. However, the resulting control is less responsive to extreme price events and does not exploit favourable or unfavourable price conditions as effectively as the DPP-based strategy. In contrast, the TOU rule is inherently tied to predefined time intervals and therefore lacks the flexibility to adapt to stochastic price fluctuations, resulting in a noticeably more rigid control behavior.
	
\begin{figure}[!t]
		\centering
		\includegraphics[width=4.in]{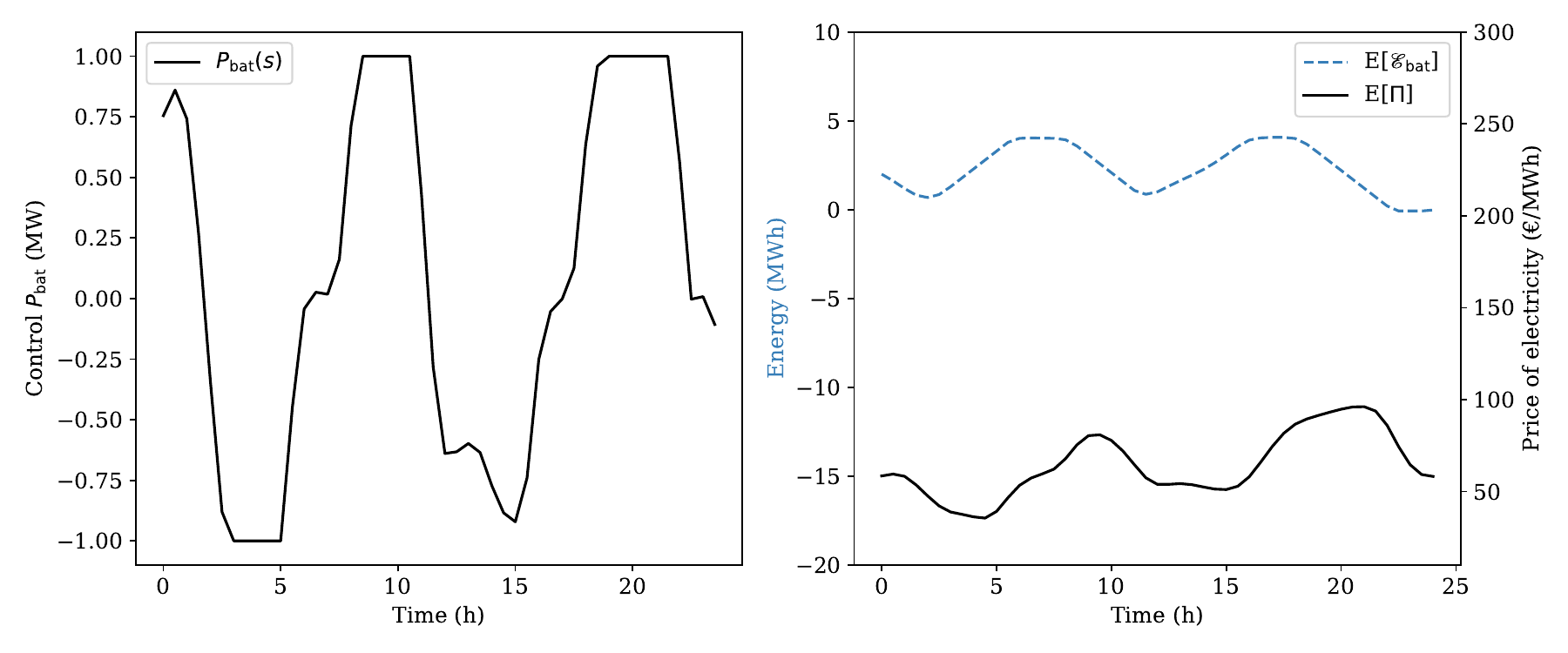}
		\caption{Trajectories of $\mathcal{E}$, $\Pi$ and $P_{\text{bat}}$ via stochastic MPC.}
		\label{fig:exp_MPC}
	\end{figure}

	In the case of stochastic MPC, Fig.~\ref{fig:exp_MPC} depicts the evolution of the battery energy $\mathcal{E}_{\text{bat}}$, the electricity price $\Pi$, and the corresponding control trajectory. As previously discussed, the state trajectory is constructed by propagating, at each MPC stage, the expected value of the state after applying only the first control action of the optimized sequence. As a consequence, the imposed constraint-handling strategy guarantees that the expectation of the state remains within the prescribed limits throughout the time horizon, ensuring the feasibility of the obtained trajectory. In addition, the resulting energy trajectory exhibits a tighter saturation of the state constraints, a behavior similar to the reduced-dimensional 1D DPP formulation. This effect is further reflected in the control profile, which frequently switches between its admissible bounds. Overall, the stochastic MPC solution shows a qualitatively consistent adaptation to price fluctuations, capturing the main economic signals in a manner comparable to the DPP-based approaches.	

	We compare now the performance of the three cases based on the proposed DPP methodology by evaluating the objective function $J$ and computing the expected value of the value function at the initial time, $v_0 = v(0,(Z_0,\Pi_0,\mathcal{E}_0))$. We also report the required computing time in minutes. This comparison is performed for several mesh configurations in each case, and the results are summarized in Table~\ref{tab:tab1}.
	
	As expected, reducing the spatial dimension of the optimality system leads to a substantial reduction in computational cost, with the 3D case proving particularly demanding. This difference is further evidenced by the cost of solving the two PDEs involved in each iteration of the final algorithm, namely the FP and HJB equations. For instance, with a structured mesh defined by $\mathbf{N} = (21,66,71)$ spatial nodes in each direction and a time step of $\Delta s = 0.5$ h, the average solution times per iteration in the 3D case are approximately $14.0$ minutes for the FP equation and $25.5$ minutes for the HJB equation. In contrast, the 2D case requires only $1.9$ seconds and $7.8$ seconds, respectively, while in the 1D case the solution times drop further to $0.19$ and $0.32$ seconds.
	
	Regarding the objective function, Table~\ref{tab:tab1} shows that the 1D case consistently yields lower profit for the PV plant, reflecting the limitations in battery energy management when the state space is overly simplified. The 2D case achieves slightly higher gains than the 3D model, although the improvement is modest. This explains the close similarity observed in the expected control actions and battery energy trajectories between the 2D and 3D formulations. 
	
	From the perspective of real-world applications, the dimension reduction technique becomes essential for rendering the approach computationally viable. Solution times in both 2D and 1D cases remain manageable with moderately fine meshes, which is particularly relevant when the problem must be solved in advance to participate in day-ahead market auctions, with the 2D formulation providing the best balance between performance and cost. Furthermore, although only expectations are reported here, it should be emphasized that the DPP also delivers a closed-loop feedback control, enabling real-time management of the PV--battery system under uncertainty. This feedback is a powerful tool to steer the system state in real time when deviations from expected values occur. Consequently, a balance must be struck between computational effort and mesh refinement, ensuring sufficiently fast resolution while retaining adequate richness in the control.

	\begin{table}
		\begin{center}
			\caption{Objective function evaluation, expectation of value function at initial time and computing times in minutes for different time and space dicretizations.}
			\label{tab:tab1}
			\begin{tabular}{llll}
				\hline
				& 3D problem & 1D Problem & 2D Problem  
				\\ 
				\hline
				$\begin{aligned}
					\Delta s &= 1 \\
					N_1 &= 21 \\
					N_2 &= 66 \\
					N_3 &=36
				\end{aligned}$
				& 
				$\begin{aligned}
					&J = 2738.1 \\
					&\mathrm{E}[v_0] = 2789.9 \\
					&\text{Time} = 653.6
				\end{aligned}$
				&
				$\begin{aligned}
					&J = 2524.1 \\
					&\mathrm{E}[v_0] = 2295.3 \\
					&\text{Time} = 2.1
				\end{aligned}$
				&
				$\begin{aligned}
					&J = 2794.3 \\
					&\mathrm{E}[v_0] = 2794.4 \\
					&\text{Time} = 2.6
				\end{aligned}$
				\\ 
				\hline
				$\begin{aligned}
					\Delta s &= 0.5 \\
					N_1 &= 21 \\
					N_2 &= 66 \\
					N_3 &=71
				\end{aligned}$
				& 
				$\begin{aligned}
					&J = 2686.4 \\
					&\mathrm{E}[v_0] = 2698.1 \\
					&\text{Time} = 2967.2
				\end{aligned}$
				&
				$\begin{aligned}
					&J = 2492.9 \\
					&\mathrm{E}[v_0] = 2288.2 \\
					&\text{Time} = 2.0
				\end{aligned}$
				&
				$\begin{aligned}
					&J = 2744.6 \\
					&\mathrm{E}[v_0] = 2682.6 \\
					&\text{Time} = 14.0
				\end{aligned}$
				\\
				\hline
				$\begin{aligned}
					\Delta s &= 0.1 \\
					N_1 &= 21 \\
					N_2 &= 131 \\
					N_3 &= 141
				\end{aligned}$
				& 
				
				&
				$\begin{aligned}
					&J = 2478.3 \\
					&\mathrm{E}[v_0] = 2311.0 \\
					&\text{Time} = 24.5
				\end{aligned}$
				&
				$\begin{aligned}
					&J = 2633.2 \\
					&\mathrm{E}[v_0] = 2554.0 \\
					&\text{Time} = 3340.9
				\end{aligned}$
				\\
				\hline
			\end{tabular}
		\end{center}
	\end{table}

	Table~\ref{tab:tab2} presents objective function evaluation and computing time in minutes for the benchmark methods. Compared with the proposed DPP-based framework, both rule-based approaches achieve performance levels similar to the 1D case. However, it should be stressed that the 1D formulation relies solely on the battery's energy content, whereas the price threshold and TOU rules also incorporate electricity prices. In contrast, when compared with the 2D and 3D cases, where richer control actions are available, our methodology consistently outperforms the rule-based methods. In terms of computing time, the 1D and 2D cases remain competitive with the benchmarks under reasonably fine meshes, underscoring the effectiveness of the dimension reduction technique. Overall, the proposed framework matches rule-based methods when limited to the simplest controls, but clearly surpasses them when higher-dimensional control strategies are considered. This underlines the strength of DPP, which delivers the optimal control at every state point, in contrast to rule-based heuristics that can only provide suboptimal decisions. 
	
	On the other hand, stochastic MPC achieves objective function values that slightly improve upon those obtained with the 1D DPP formulation, but remain significantly below those of the 2D and 3D cases. This outcome is consistent with the qualitative similarities observed between the MPC and 1D DPP control behaviors. In contrast, the computational burden of stochastic MPC is notably higher: its solution times exceed those of the 1D and even some of the 2D DPP formulations, further highlighting the importance of the proposed dimension-reduction strategy within the DPP framework. Moreover, it is important to note that the MPC scheme yields an open-loop control trajectory that is sequentially updated in real time, rather than a feedback control law defined over the state space. As a consequence, the resulting control lacks the structural richness of the DPP-based solutions and is inherently less suited for day-ahead market participation, where decisions must be committed in advance.
	\begin{table}
		\begin{center}
			\caption{Objective function evaluation and computing times in minutes for benchmark methods.}
			\label{tab:tab2}
			\begin{tabular}{lll}
				\hline
				Price thresholds & TOU & Stochastic MPC
				\\ 
				\hline
				$\begin{aligned}
					&J = 2363.9 \\
					&\text{Time} = 7.9
				\end{aligned}$
				&
				$\begin{aligned}
					&J = 2277.8 \\
					&\text{Time} = 7.5
				\end{aligned}$
				&
				$\begin{aligned}
					&J = 2531.8 \\
					&\text{Time} = 67.9
				\end{aligned}$
				\\
				\hline
			\end{tabular}
		\end{center}
	\end{table}
	
	In addition to the previous results, we analyze the sensitivity of Algorithm~\ref{alg:alg1} to the choice of the initial penalty parameter $\lambda^0$. Within the proposed penalty-multiplier framework, this parameter governs the strength of constraint penalization during the initial iterations, when the algorithm effectively operates in a pure penalty regime. Table~\ref{tab:tab3} reports, for the 1D, 2D, and 3D formulations, the number of iterations required to reach convergence together with the relative $L^2$-difference between control solutions obtained for consecutive values of $\lambda^0$. For each spatial dimension, the remaining algorithmic parameters are fixed, while $\lambda^0$ is varied over several orders of magnitude. In particular, the sequence of values for $\lambda^0$ is chosen to decrease progressively towards a prescribed lower bound $\lambda^*$, which defines the minimum admissible value of the penalty parameter within the algorithm. This lower bound is introduced to prevent ill-conditioning issues observed for excessively small values of $\lambda$, which may impede convergence in practice. Although the specific value of $\lambda^*$ depends on the problem dimension, it remains of the same order of magnitude across the 1D, 2D, and 3D cases.
	
	The results show a consistent reduction in the number of iterations as $\lambda^0$ decreases, indicating faster convergence when starting from a stronger penalization regime. However, as $\lambda^0$ approaches $\lambda^*$, this reduction becomes less pronounced, and the number of iterations tends to stabilize. From the accuracy standpoint, the computed solutions remain highly stable across all tested configurations. In particular, the relative differences in the control between successive configurations are negligible. Numerical results also showed that the corresponding objective function values coincide up to the fifth significant digit. This confirms that the choice of $\lambda^0$ primarily affects convergence speed rather than solution quality. These findings are particularly relevant in higher-dimensional settings as the 3D case, where each iteration entails a significant computational cost due to the solution of coupled FP-HJB equations. Overall, these results support the robustness of the proposed penalty-multiplier strategy with respect to the initial selection of the penalty parameter.

	\begin{table}
	\begin{center}
		\caption{Impact of the initial penalty parameter $\lambda^0$ on convergence and control stability across different problem dimensions.}
		\label{tab:tab3}
		\begin{tabular}{cccc}
			\hline
			Spatial dimension & $\lambda^0$ & Iterations & \shortstack{Relative $L^2$-difference in control \\ (with respect to previous $\lambda^0$)}
			\\
			\hline
			3D 
			& 
			${\setlength{\jot}{1pt}\begin{aligned}
				10^3& \\ 10^2& \\ 50& \\ 20& \\ 8& \\ 6& 
			\end{aligned}}$
			& 
			${\setlength{\jot}{1pt}\begin{aligned}
				57& \\ 44& \\ 42& \\ 39& \\ 38& \\ 39& 
			\end{aligned}}$ 
			& 
			${\setlength{\jot}{1pt}\begin{aligned}
					-& \\ 0.040& \\ 0.036& \\ 0.026& \\ 0.035& \\ 0.019& 
			\end{aligned}}$
			\\
			\hline
			1D 
			& 
			${\setlength{\jot}{1pt}\begin{aligned}
					10^3& \\ 10^2& \\ 50& \\ 20& \\ 8& \\ 6& 
			\end{aligned}}$
			& 
			${\setlength{\jot}{1pt}\begin{aligned}
					358& \\ 348& \\ 345& \\ 345& \\ 211& \\ 211& 
			\end{aligned}}$ 
			& 
			${\setlength{\jot}{1pt}\begin{aligned}
					-& \\ 0.221& \\ 0.037& \\ 0.051& \\ 0.621& \\ 0.023& 
			\end{aligned}}$
			\\
			\hline
			2D 
			& 
			${\setlength{\jot}{1pt}\begin{aligned}
					10^3& \\ 10^2& \\ 50& \\ 20& \\ 8& \\ 6& 
			\end{aligned}}$
			& 
			${\setlength{\jot}{1pt}\begin{aligned}
					62& \\ 52& \\ 50& \\ 46& \\ 45& \\ 44& 
			\end{aligned}}$ 
			& 
			${\setlength{\jot}{1pt}\begin{aligned}
					-& \\ 0.051& \\ 0.041& \\ 0.066& \\ 0.051& \\ 0.051& 
			\end{aligned}}$
			\\
			\hline
		\end{tabular}
	\end{center}
	\end{table}

	\section{Conclusion} \label{sect:VI_conclusion}
	\noindent This paper builds on the framework introduced in \cite{BP} for wind farms, extending it to the joint optimization of real-time energy management and day-ahead market participation for PV plants with storage. The consistency of the results across these two technologically different settings highlights the robustness and generality of the proposed continuous-time stochastic control framework. Key to our approach is the reformulation of the original stochastic optimal control problem as a deterministic one via the FP equation. The resulting optimality system is fully coupled, consisting of a forward FP equation and a backward HJB equation. State constraints are enforced through a Lagrange multiplier that appears naturally in the HJB equation. We then design and implement a specialized numerical algorithm to solve this coupled FP-HJB system efficiently.
	
	While the previous study in \cite{BP} focused on demonstrating the feasibility and effectiveness of the continuous-time Dynamic Programming approach, the present work advances the framework by explicitly addressing its computational scalability. In particular, we introduce and validate a dimension reduction strategy that preserves the structure of the original optimality system while enabling a systematic trade-off between control expressiveness and computational effort. We decompose the state vector into two parts, isolating a subset of uncontrollable variables whose marginal distribution can be computed independently. The remaining (controllable) variables then define a reduced dimension stochastic optimal control problem.
	
	Numerical experiments on a realistic PV-plant case study show that reduced-dimensional formulations retain the essential qualitative features of the optimal policy and achieve competitive objective function values, while drastically reducing computation times. This marks a significant step forward relative to earlier work, bringing the continuous-time DPP framework closer to practical real-time deployment and market-oriented applications. Comparative benchmarking against two rule-based methods, namely, price threshold and TOU, as well as stochastic MPC further highlights the effectiveness of the proposed approach. In particular, the reduced-dimensional formulations achieve a favorable compromise between computational efficiency and control performance, yielding solutions that are competitive with, and in some cases superior to, the benchmark methods. Unlike rule-based approaches, the proposed methodology preserves an optimal control structure, and unlike stochastic MPC, it provides a feedback control law defined over the state space. As a result, the reduced DPP formulations are well suited for real-time operation while simultaneously retaining the ability to support day-ahead market participation.

	Despite the proposed dimension-reduction strategy, the methodology remains fundamentally constrained by the curse of dimensionality as the number of state variables increases. In particular, when a large set of states must be directly controlled, the computational effort required to solve the associated partial differential equations can become prohibitive. Thus, although reduced-dimensional formulations can still generate optimal feedback policies suitable for both real-time operation and day-ahead market participation, further increases in model complexity may make simpler approaches more attractive from a computational standpoint. In this sense, while the proposed framework offers a comprehensive and theoretically optimal solution, it may be outperformed in terms of computing time by less expressive methods--such as discrete-time optimization strategies--when applied to very high-dimensional systems.

	Future research will primarily focus on addressing the computational challenges inherent to high-dimensional stochastic control problems. In this direction, data-driven and machine-learning techniques offer promising complementary tools for approximating value functions or feedback control laws, while preserving the continuous-time Dynamic Programming structure of the framework. Parallel developments aim to enrich the problem formulation to better capture real-world operating conditions. These include incorporating additional sources of uncertainty and extending the framework to multi-agent or networked energy systems. Such enhancements would enable the methodology to more accurately reflect practical energy-market environments, where uncertainty is multifaceted and decision-making is often interdependent across multiple assets or market participants.

	\section*{Declaration of competing interest}
	The authors declare that they have no known competing financial interests or personal relationships that could have appeared to influence the work reported in this paper.

	\section*{Acknowledgement}
	This work has been partially supported by FEDER, Ministerio de Ciencia e Innovaci\'on-AEI research project PID2021-122625OB-I00.

\end{document}